\documentclass[12pt,a4paper,english]{article}
\usepackage{graphicx} 
\usepackage[T1]{fontenc}
\usepackage{amssymb}
\makeatletter

\pdfpageheight\paperheight
\pdfpagewidth\paperwidth

\usepackage{german}
\usepackage[ansinew]{inputenc}
\usepackage{amsmath}
\usepackage{amsfonts}
\usepackage[all]{xy} %f"ur kommutative Diagramme, xypics 
\catcode`\"=\active \let"=\" \let\3=\ss
\def\heute{\number\day.\space\ifcase\month\or Januar\or Februar\or M"arz\or April
\or Mai\or Juni\or Juli\or August\or September\or Oktober\or November
\or Dezember\fi\space\number\year}
\newcommand\exa{\nopagebreak \begin{flushleft}\smallskip\nopagebreak
                \begin{minipage}[t]{8cm}\sloppy}
\newcommand\exb{\end{minipage}\kern 1cm\begin{minipage}[t]{8cm}\sloppy }
\newcommand\exc{\end{minipage}\kern -3cm \smallskip\end{flushleft}}

\makeatother

\usepackage{babel}
\begin{document}
\date {}

\title{An Algorithm to Find Rational Points on Elliptic Curves Related to the Concordant Form Problem}
\author{Hagen Knaf, Erich Selder, Karlheinz Spindler}

\maketitle
{\bf Abstract} We derive an efficient algorithm to find solutions to Euler's 
concordant form problem and rational points on elliptic curves associated with 
this problem. \\ 
 
{\bf Keywords} algorithmic number theory, elliptic curves, concordant forms \\

{\bf AMS Subject Classification} 11Y16, 11G05, 14H52  
$$\phantom{a}$$ 

\section{Introduction}

It is well known that the determination of the Mordell-Weil group of an elliptic curve is a 
difficult problem. Apart from the torsion subgroup, which can be calculated rather easily 
using the Lutz-Nagell-Theorem ([11], [16], cf. [20], [21]) and for which very good general 
information is given by Mazur's Theorem (cf. [12]), obtaining information on the rational 
points of a rationally defined elliptic curve is hard. Even if the elliptic curve is 
explicitly given, both the calculation of the rank of the Mordell-Weil group 
and the determination of explicit solutions (generators of this group) are 
difficult problems. A famous example due to Zagier (cf. [25]; also see [9], p.5, 
Fig. 1.3), which illustrates the problems, is the task of explicitly showing that 
$n=157$ is a congruent number by determining the sides of a rational triangle with 
area $157$, which consist of fractions with more than $25$ decimal places in both 
the numerator and the denominator. As a byproduct of these calculations one 
easily obtains nontrivial rational solutions of the equation $y^{2}=x(x-157)(x+157)$,
the smallest of which is given by 
\def\spc{\hskip -1.5 true mm}
\begin{eqnarray*}
x &\spc =\spc& -166136231668185267540804/2825630694251145858025,\\
y &\spc =\spc& 167661624456834335404812111469782006/150201095200135518108761470235125.
\end{eqnarray*}
An even more exciting example was discovered by Bremner and Cassels in [1] who 
considered the family of curves with equations $y^{2}=x(x^{2}+p)$ with prime 
numbers $p\equiv5$ (mod $8$). In the introduction to this paper the authors mention 
an observation of S. Lang which asserts that in the tables of elliptic curves and 
generators of the Mordell-Weil group known at that time (around 1982) the numerator 
and the denominator of the $x$-value of smallest nontrivial solutions to elliptic 
curve equations with integral coefficients are never much greater than the square 
of the discriminant. Note that the above example of Zagier exemplifies this 
observation, since the discriminant of the example is $\triangle=16\cdot157^{6} = 
239617149303184$. \smallskip 

The most spectacular example calculated in [1] (see also [20], Chap. X, §6, Remark 6.3) 
is given by the prime number $p=877$. The generator of the infinite part of the 
Mordell-Weil group in this example is given by 
$$x = \frac{375494528127162193105504069942092792346201}{62159877768644257535639389356838044100}$$ 
and 
$$y = \frac{256256267988926809388776834045513089648669153204356603464786949}{
   490078023219787588959802933995928925096061616470779979261000}$$ 
so that both the numerator and the denominator of $x$ are greater than
$\Delta^{4}$, thus providing a clear counterexample to Lang's observation.\smallskip 

In this paper we present an algorithm which enables us to calculate
a large number of examples of smallest solutions to certain classes 
of elliptic curve equations. Many of the explicitly computed examples
belonging to classes of elliptic curves other than those of the paper 
of Bremner and Cassels (namely those corresponding to Euler's concordant 
forms) also contradict Lang's assertion. We are not interested in general 
discussions on the determination of the rank of the Mordell-Weil group, but focus 
on finding explicit solutions. So in general we assume the given
elliptic curves to have positive Mordell-Weil-rank. The method for
finding solutions is a simple search loop of quadratic complexity,
but the key point is the reduction of the equations to be considered
to simpler ones by some technical tricks. \smallskip 

The first (and most important) reduction is a well-known descent
procedure in which the original elliptic curve $E$ is substituted
by a homogeneous space $Q$ over $E$. The equations defining this
homogeneous space can be simplified by Newton's method of parametrizing
a quadric by means of a projection from a fixed point onto that quadric.
These two concepts already enable us to formulate a simple algorithm,
allowing us to find explicit solutions for a large class of curves.
\smallskip 

The core of our approach is an improvement of the above algorithm which 
makes it possible to obtain solutions of twice the complexity as compared 
to the simpler algorithm. This improvement is achieved by invoking a technical 
condition which is satisfied in many situations.\smallskip 

The potential of this improved algorithm will be illustrated by some series of 
examples. For example, this algorithm can be used to find the smallest nontrivial 
integers $A,B,C,D$ satisfying the two equations $A^2-373B^2 = C^2$ and 
$A^2+373B^2=D^2$, which are found to be 

\def\spc{\hskip 0 true mm}
{\footnotesize \begin{eqnarray*}
A&\spc =\spc&6464736286838262275566375140640125524476830394378258160144359151221846588162921 \\  
B&\spc =\spc&214402886988423616335778394508029972671920911384749815755228436417174376951980 \\ 
C&\spc =\spc&4964526988887992094202607668810309975770378526931158358479760499172740751760929 \\ 
D&\spc =\spc&7677180621382399924131415436519959747090354653821331133153517438341892919535729 
\end{eqnarray*} } 

\noindent 
and have up to $79$ decimals. (The existence of such a solution shows explicitly that 
$373$ is a congruent number.) Based on an analysis of the examples to be presented 
later in the paper, it may be possible to gain some insight into the behaviour of the 
solutions depending on the input parameters of the given elliptic curves, but many 
open questions remain. We see our contribution mainly in providing a tool which allows 
one to systematically find rational points on certain classes of elliptic curves. 
\smallskip 

The paper is organized as follows. After this introduction we present some preparatory 
material on elliptic curves and quadratic forms. The most important part of this section 
will be the description of the abstract descent procedure and its concrete realization 
for our purposes. In a subsequent section we illustrate this descent procedure by some examples, 
some of which will be taken up again in the series of examples discussed at the end of the 
paper. Subsequently, we describe our algorithm in detail. In particular, we point out the 
critical choices to be made during the execution of the algorithm and the conditions which 
are necessary to make it work. A further section contains some series of examples showing
the power of the algorithm. Finally, we close with a critical discussion of the algorithm, 
examine consequences of the examples generated by the algorithm and give an outlook to 
further questions to be pursued.

\section{General background}

\subsection{Elliptic curves}

\subsubsection{General notions}

Throughout this paper we consider most of our geometric objects to
be defined over the field $\mathbb{Q}$ of rational numbers. Sometimes
we argue geometrically and regard the objects to be defined over an 
algebraic closure $\overline{\mathbb{Q}}$ or over some number field. 
\smallskip 

An elliptic curve defined over $\mathbb{Q}$ is a plane projective
curve $E$ whose affine part is given by a Weierstra\ss{} equation of the
form $y^{2}=P(x)$ where $P\in\mathbb{Q}[x]$ is a polynomial of degree $3$ 
with nonzero discriminant. Note that such a curve always has a 
single smooth point at infinity which is defined over $\mathbb{Q}$. The 
geometric points of $E$ carry the structure of an algebraic group, the 
group law being given by Newton's well-known secant and tangent construction. 
We fix the structure such that the point at infinity is the neutral
element of this group law. The group law restricts to the set $E(\mathbb{Q})$
of rational points on $E$, the so-called Mordell-Weil group of the
elliptic curve $E$. By the important theorem of Mordell and Weil (cf. [13],
[22]) this group is a finitely generated abelian group, hence is isomorphic
(as an abstract group) to $\mathbb{Z}^{r}\times T$ where $r=\hbox{rank} 
(E(\mathbb{Q}))$ is the rank of $E(\mathbb{Q})$ and $T$ is a torsion group, 
i. e., a finite abelian group. \smallskip 

If $E$ is explicitly given then it is relatively easy to compute
the torsion part of the Mordell-Weil group by the theorem of Lutz-Nagell. 
Moreover, the famous theorem of Mazur tells us that there are only very 
few possibilities for a finite abelian group to occur as the 
torsion subgroup of the Mordell-Weil group of an elliptic curve. In fact, $T$ is 
either isomorphic to $\mathbb{Z}/n\mathbb{Z}$ where $1\leq n\leq12$ , $n\neq11$, 
or else is isomorphic to $\mathbb{Z}/2\mathbb{Z}\times\mathbb{Z}/2n\mathbb{Z}$
where $1\leq n\leq 4$. \smallskip  

In contrast to the torsion part, the determination of the rank of the
Mordell-Weil group is a much deeper problem. Moreover, it is very
difficult to find non-torsion solutions even if one knows in advance
that the rank is positive, i. e., that there exist infinitely many rational
points on $E$. This is due to the fact that even if the coefficients
of the elliptic curve are small (measured for example in terms of a naive
notion of height in projective space) the smallest non-torsion point on $E$ 
may be very large. 

\subsubsection{Elliptic curves corresponding to concordant forms}

In the sequel we will deal with elliptic curves $E$ defined over
the rationals with torsion subgroup containing $\mathbb{Z}/2\mathbb{Z}
\times\mathbb{Z}/2\mathbb{Z}$. This is equivalent to saying that $E$ is 
defined by an affine equation of the form $y^2 = (x-e_1)(x-e_2)(x-e_3)$ 
with pairwise different rational numbers $e_1,e_2,e_3\in\mathbb{Q}$.
We denote this elliptic curve by $E_{e_1,e_2,e_3}$. By simple
rational transformations we can clear denominators and translate the
$x$-coordinates such that the equation becomes $y^2=x(x+M)(x+N)$ 
with different nonzero integers $M,N\in\mathbb{Z}$; we denote this curve 
by $E_{M,N}$. (In projective notation we use homogeneous coordinates 
$(T,X,Y)$ and write the points of the affine part of $E_{M,N}$ as 
$(x,y) = (X/T, Y/T)$.) In this form the elliptic curve correponds to Euler's
concordant form problem (cf. [3], [17]), which is the problem of finding 
nontrivial integral solutions of the system of two quadratic equations 
$X_0^2+MX_1^2=X_2^2$ and $X_0^2+NX_1^2=X_3^2$. (See [2], Ch. 8, for the 
representation of elliptic curves as intersections of quadrics.) The intersection 
of these two quadrics in projective three-space $\mathbb{P}^{3}(\overline{
\mathbb{Q}})$ will be denoted by $Q_{M,N}$. For future reference, let us 
fix the notations\\ 
\parbox{.03\textwidth}{(1)}
\parbox{.97\textwidth}{
\begin{align*}  
  E_{M,N} & = \{ (T:X:Y)\in\mathbb{P}^2(\mathbb{\overline{Q}})\mid 
       TY^2=X(X\! +\! TM)(X\! +\! TN)\}, \\
  Q_{M,N} & = \{ (X_0:X_1:X_2:X_3)\in\mathbb{P}^3(\mathbb{\overline{Q}})\mid 
       X_0^2\! +\! MX_1^2=X_2^2,\ X_0^2\! +\! NX_1^2=X_3^2\}. 
  \end{align*} }\\ 
We note that $Q_{M,N}$ defines an abstract elliptic curve (i.e., a smooth projective 
curve of arithmetic genus one) which is isomorphic to the curve $E_{M,N}$. In fact, 
the mappings $F:Q_{M,N}\rightarrow E_{M,N}$ given by 
$$F\left(\begin{array}{c} X_0\\ X_1\\ X_2\\ X_3\end{array}\right) = 
  \left(\begin{array}{c} NX_2-MX_3+(M-N)X_0\\ MN(X_3-X_2)\\ MN(M-N)X_1
  \end{array}\right) \leqno{(2)}$$
and $G:E_{M,N}\rightarrow Q_{M,N}$ given by 
$$G\left(\begin{array}{c} T\\ X\\ Y\end{array}\right) = \left(\begin{array}{c}
  -(X+MT)(Y^2-M(X+NT)^2)\\ 2Y(X+NT)(X+MT)\\ -(X+MT)(Y^2+M(X+NT)^2)\\
  -(X+NT)(Y^2+N(X+MT)^2) \end{array}\right)\leqno{(3)}$$ 
extend to completely defined biregular mappings which are mutually
inverse; cf. [19]. Note that for the equation of $E_{M,N}$ we may assume that 
$M>0$ and $N<0$ after a trivial change of 
variables. Furthermore, we may write $M=pk$ and $-N=qk$ with coprime natural 
numbers $p,q\in\mathbb{N}$ and a squarefree natural number $k\in\mathbb{N}$. 
In this form the equation is closely related to the $\theta$-congruent number 
problem in the sense of Fujiwara (cf. [4]; see also [7], [24]). \smallskip  

Note that the two-torsion points $(0,0)$, $(-M,0)$, $(-N,0)$ together
with the point at infinity on $E_{M,N}$ correspond to the trivial
solutions $(1,0,\pm1,\pm1)$ of the concordant form equations and
to the trivial solution (the degenerated triangle) of the $\theta$-congruent
number problem. Let us point out that by Mazur's Theorem the only further
torsion points could be either $4$- or $8$-torsion points or else
$3$- or $6$-torsion points. Moreover, these additional torsion points
can occur only in very rare cases. In fact, in the form $E_{pk,-qk}$
of the elliptic curve such torsion points can only occur when $k=1,2,3$
or $6$ (cf. [5], [19]). Since the torsion points can be calculated easily 
we will ignore them in the subsequent considerations in which we develop methods 
for finding non-torsion points on elliptic curves of the considered form. 
\smallskip 

\subsection{Quadratic Forms}

In the following sections we frequently consider ternary quadratic
forms, in most cases of a rather special form. We will collect some
facts about these forms needed later. The quadratic forms which occur 
in the following considerations are of the form
$F(X_0,X_1,X_2)=a_{00}X_0^{2}+a_{01}X_0X_1+a_{11}X_1^2+a_{22}X_2^2$
with coprime integer coefficients $a_{ij}\in\mathbb{Z}$. In most cases 
the forms are already in diagonal form (i.e., $a_{01}=0$) or are 
transformed to diagonal form in the course of our investigations. 

\subsubsection{Solvability criterion}

Let $F(X_0,X_1,X_2)=a_{00}X_0^2+a_{11}X_1^2+a_{22}X_2^2$ be a ternary 
diagonal form with integer coefficients. Then by standard techniques we 
can transform this equation into one for which the product $a_{00}
a_{11}a_{22}$ is squarefree (equivalently, such that the three coefficients 
are squarefree and pairwise coprime). The following criterion for the
existence of a nonzero integer solution $(x_0,x_1,x_2)\in\mathbb{Z}^{3}$ 
dates back to Legendre (cf. [10]). \bigskip

{\bf Legendre Criterion:} The equation $a_{00}X_0^2+a_{11}X_1^2+a_{22}X_2^2
=0$ has a nonzero solution if and only if not all the coefficients have the 
same sign (which is equivalent to saying that there is a solution over the 
real numbers) and, in addition, for all permutations $(i,j,k)$ of $(0,1,2)$ 
the number $-a_{ii}a_{jj}$ is a quadratic residue modulo $a_{kk}$.\bigskip 

With this criterion one can easily check whether or not a given ternary 
quadratic form has a nontrivial solution. In addition, by the work of Holzer 
and his followers (cf. [6], [14], [23]) there are algorithms to find explicit 
solutions which terminate with a complexity which is known {\it a priori}  
depending on the coefficients. (Incidentally, Holzer's Theorem occurs both 
as a tool within our algorithm and as a model for our overall approach, namely, 
to produce explicit solutions to certain types of equations in which the mere 
existence of solutions is known by other, more abstract, methods.) \bigskip  

{\bf Holzer's Theorem:} If $F(X_0,X_1,X_2)=a_{00}X_0^2+a_{11}X_1^2+a_{22}X_2^2$
is a quadratic form with pairwise coprime and squarefree coefficients such 
that the equation $F(X_0,X_1,X_2)=0$ has a nontrivial solution, then
there exists such a solution $(x_0,x_1,x_2)$ satisfying the inequalities
$|x_0|<\sqrt{|a_{11}a_{22}|}$ , $|x_1|<\sqrt{|a_{00}a_{22}|}$ and 
$|x_{2}|<\sqrt{|a_{00}a_{11}|}$.\bigskip  

\subsubsection{Parametrization of quadratic forms}

In the sequel we will make use of the fact that there is a systematic
approach to finding all solutions of a quadratic form, provided one 
fixed solution is known. The technique is a well-known construction,
dating back to Newton, which geometrically uses the fact that any
line through the fixed point has a uniquely determined second point
of intersection with the quadric, whose coordinates are rational 
expressions in the coefficients of the given quadric and the coordinates 
of the fixed point, where the slope of the line serves as a parameter. We 
summarize the calculations, writing down the parametrization using projective 
coordinates.\bigskip

{\bf Lemma:} Consider the projective quadric $Q=\{(x_0,x_1,x_2)\in\mathbb{P}^{2}
\,|\,F(X_0,X_1,X_2)$ $=0\}$ where $F(X_0,X_1,X_2)=a_{00}X_0^2+a_{01}X_0X_1+a_{11}X_1^2
+a_{22}X_2^2$. The set of all points $(x_0,x_1,x_2)\in Q$ with $x_2\not= 0$ can be 
parametrized by the rational mapping $\varPhi:\mathbb{P}^1\rightarrow Q$ given by 
$\varPhi(\xi_0,\xi_1) = (\varphi_0(\xi_0,\xi_1),\varphi_1(\xi_0,\xi_1),\varphi_2(\xi_0,\xi_1))$ 
where 

\parbox{.05\textwidth}{(4a)}
\parbox{.95\textwidth}{
  \begin{alignat*}{3}
  & X_0\ &&=\ \varphi_0(\xi_0,\xi_1)\ &&=
       \ a_{11}x_0\xi_0^2 - 2a_{11}x_1\xi_0\xi_1 - (a_{01}x_1+a_{00}x_0)\xi_1^2,  \\
  & X_1\ &&=\ \varphi_1(\xi_0,\xi_1)\ &&=  
       \ (-a_{11}x_1-a_{01}x_0)\xi_0^2 - 2a_{00}x_0\xi_0\xi_1 + a_{00}x_1\xi_1^2, \\ 
  & X_2\ &&=\ \varphi_2(\xi_0,\xi_1)\ &&=  
       \ a_{11}x_2\xi_0^2 + a_{01}x_2\xi_0\xi_1 + a_{00}x_2\xi_1^2.  
  \end{alignat*} } \\
whereas if $x_2=0$ and $x_1\not= 0$ the following parametrization can be used: \\
\parbox{.05\textwidth}{(4b)}
\parbox{.95\textwidth}{
  \begin{alignat*}{3}
  & X_0\ &&=\ \varphi_0(\xi_0,\xi_1)\ &&=
       \ -(a_{01}x_1+a_{00}x_0)\xi_0^2+a_{22}x_0\xi_1^2,  \\
  & X_1\ &&=\ \varphi_1(\xi_0,\xi_1)\ &&=  
       \ a_{00}x_1\xi_0^2 + a_{22}x_1\xi_1^2, \\ 
  & X_2\ &&=\ \varphi_2(\xi_0,\xi_1)\ &&=  
       \ -(a_{01}x_1+2a_{00}x_0) \xi_0\xi_1.  
  \end{alignat*} } \\
The proof is an elementary and simple calculation, which we omit.
\bigskip 

{\bf Remark:} We did not specify a base field over which the projective
spaces and all the coordinates are defined, since the lemma is valid for any 
field. We will apply this lemma in the situation that the quadric is defined 
over the rational numbers; i.e., the coefficients $a_{00},a_{01},a_{11},a_{22}$ 
as well as the coordinates of the fixed point $(x_0,x_1,x_2)$ and the parameters 
$(\xi_0,\xi_1)$ are rational numbers. In this situation we can always assume that 
all these data are, in fact, integers such that the coefficients $a_{00}, a_{01}, 
a_{11}, a_{22}$ are coprime, the coordinates $x_0,x_1,x_2$ are coprime and the 
parameters $\xi_0,\xi_1$ are coprime.\\

\subsubsection{Pairs of quadrics with separated variables}

In the situation of the following discussion we will sometimes be
given two quadrics in projective three-space of a special form in which 
the variables of the quadrics will be separated in some sense. More precisely,
the quadrics $Q_{1}$ and $Q_{2}$ will have equations of the form \\
\parbox{.05\textwidth}{(5)}
\parbox{.95\textwidth}{
  \begin{alignat*}{3}
    & Q_1: && \quad a_{00}X_0^2 + a_{01}X_0X_1 + a_{11}X_1^2 + a_{22}X_2^2\ &&=\ 0, \\
    & Q_2: && \quad b_{00}X_0^2 + b_{01}X_0X_1 + b_{11}X_1^2 + b_{33}X_3^2\ &&=\ 0
  \end{alignat*} } \\
so that the equations for the quadrics share two of the four variables such that 
the other two variables occur only as squares. We always assume the quadrics to have 
rational solutions, so by means of a fixed point $(x_0,x_1,x_2)$ on $Q_1$ we can 
parametrize the points on the first quadric $Q_1$, viewed as a quadric in the projective 
plane with coordinates $(X_0,X_1,X_2)$, and we can substitute the parametrizations for 
the common variables $X_0,X_1$ into the second quadric, thus yielding a function of the 
variables $\xi_0,\xi_1,X_3$ which is of degree $4$ in the parameter variables $\xi_0,\xi_1$ 
and of pure degree $2$ in the third variable $X_3$. Let us summarize the result of the 
calculations. \\

{\bf Lemma:} Let $Q_1,Q_2$ be two quadratic forms as above, let $(x_0,x_1,x_2)$ be a 
point on $Q_1$ and let \\ 
\parbox{.05\textwidth}{(6)}
\parbox{.95\textwidth}{
  \begin{alignat*}{3}
   & X_0\ &&=\ \varphi_0(\xi_0,\xi_1)\ &&=\ \alpha_{00}\xi_0^2 + \alpha_{01}\xi_0\xi_1 
            + \alpha_{11}\xi_1^2 \\ 
   & X_1\ &&=\ \varphi_1(\xi_0,\xi_1)\ &&=\ \beta_{00}\xi_0^2 + \beta_{01}\xi_0\xi_1 
            + \beta_{11}\xi_1^2 \\
   & X_2\ &&=\ \varphi_2(\xi_0,\xi_1)\ &&=\ \gamma_{00}\xi_0^2 + \gamma_{01}\xi_0\xi_1 
            + \gamma_{11}\xi_1^2
  \end{alignat*} } \\
be the parametrization of $Q_1$ projecting from $(x_{0},x_{1},x_{2})$.
Then the substitution of $\varphi_0$ and $\varphi_1$ into the quadric 
$Q_2$ gives the equation
$$Q_{2}:\quad B_{40}\xi_0^4 + B_{31}\xi_0^3\xi_1 + B_{22}\xi_0^2\xi_1^2 
  + B_{13}\xi_0\xi_1^3 + B_{04}\xi_1^4 + b_{33}X_3^2 \leqno{(7)}$$
where \\
\parbox{.03\textwidth}{(8)}
\parbox{.97\textwidth}{
  \begin{alignat*}{2}
  & B_{40} &&= b_{00}\alpha_{00}^{2} + b_{01}\alpha_{00}\beta_{00} + b_{11}\beta_{00}^{2} \\
  & B_{31} &&= 2b_{00}\alpha_{00}\alpha_{01} + b_{01}(\alpha_{00}\beta_{01}\! +\!\alpha_{01}
                 \beta_{00}) + 2b_{11}\beta_{00}\beta_{01} \\
  & B_{22} &&= b_{00}(2\alpha_{00}\alpha_{11}\! +\!\alpha_{01}^{2}) + b_{01}(\alpha_{00}\beta_{11}
                 \! +\!\alpha_{01}\beta_{01} + \alpha_{11}\beta_{00}) + b_{11}(2\beta_{00}\beta_{11}
                 \! +\!\beta_{01}^{2}) \\ 
  & B_{13} &&= 2b_{00}\alpha_{01}\alpha_{11} + b_{01}(\alpha_{01}\beta_{11}\! + \!\alpha_{11}
                 \beta_{01})+2b_{11}\beta_{01}\beta_{11} \\ 
  & B_{04} &&= b_{00}\alpha_{11}^{2}+b_{01}\alpha_{11}\beta_{11}+b_{11}\beta_{11}^{2}
  \end{alignat*} } \\
Again, the proof is just an easy calculation and is omitted.\\
\\
{\bf Corollary:} If the two quadrics $Q_1$ and $Q_2$ are diagonal
(i.e., if $a_{01}=b_{01}=0$) and if one of the coordinates $x_0$ or
$x_1$ of the fixed point is zero, then the substituted form of
$Q_2$ is biquadratic in $(\xi_0,\xi_1)$.\smallskip 

In fact, in this situation we have $\alpha_{00}=\alpha_{11}=\beta_{01}=0$
if $x_{0}=0$ or else $\alpha_{01}=\beta_{00}=\beta_{11}=0$ if $x_{1}=0$.
In both cases the two coefficients $B_{31}$ and $B_{13}$ vanish.

\subsection{Two-descent}

\subsubsection{General theory}

We are interested in explicitly finding rational points on the elliptic
curves $E_{e_1,e_2,e_3} = \{ (x,y)\in\mathbb{Q}^2\mid y^2=(x-e_1)(x-e_2)(x-e_3)\}$ 
as defined in 2.1.2. We may restrict our considerations to the case that the 
numbers $e_i$ are integers. According to Silverman (cf. [20], Chap. X, Remark 
3.4) the \emph{existence} of rational points may be checked by deciding whether 
certain homogeneous spaces over $E_{e_1,e_2,e_3}$ are trivial in the sense of 
[20], Chap. X, §3, Prop. 3.3. But these homogeneous spaces are also useful for 
calculating rational points explicitly. Since for our purposes we do not make 
use of the important homological criteria for the existence of solutions 
(local-global principle, Selmer group, Tate-Shafarevich group etc.), we only recall 
the principal facts on the determination of homogeneous spaces belonging to rational 
points on $E_{e_1,e_2,e_3}$. Here we may restrict our attention to the isogeny given 
by the multiplication-by-$2$-map on $E_{e_1,e_2,e_3}$. We use the notations of [20]  
(see Chap. X, Section 1 and Example 4.5.1). The same calculations can be found in [8]  
(see Chap. IV, Section 3). \\

First we observe that any element of $\mathbb{Q}^{*}/(\mathbb{Q}^{*})^{2}$ can be 
uniquely represented by a squarefree integer. Any pair $(b_1,b_2)\in\mathbb{Q}^{*}/(
\mathbb{Q}^{*})^{2}\times\mathbb{Q}^{*}/(\mathbb{Q}^{*})^{2}$ defines the intersection 
of two quadrics $Q_{e_1,e_2,e_3,b_1,b_2} = Q_1\cap Q_2$ in projective three-space by \\
\parbox{.05\textwidth}{(9)}
\parbox{.95\textwidth}{
  \begin{alignat*}{3}
  & Q_1: &&\quad b_1X_1^2 - b_2X_2^2 + (e_1-e_2)X_0^2     &&=\ 0, \\
  & Q_2: &&\quad b_1X_1^2 - b_1b_2X_3^2 + (e_1-e_3)X_0^2\ &&=\ 0.
  \end{alignat*} } \\
(Cf. [20], Chap. X, Section 1.) This intersection defines an abstract elliptic curve which 
in general does not have any rational point, but which is a twist of the given curve; i.e., 
it becomes isomorphic to $E_{e_1,e_2,e_3}$ over an algebraic closure $\overline{\mathbb{Q}}$.
This curve defines a homogeneous space over $E_{e_1,e_2,e_3}$, hence an element in the 
Weil-Chatelet group $WC(E_{e_1,e_2,e_3},\mathbb{Q})$. More precisely, since the above 
construction stems from the multiplication by $2$ on $E_{e_1,e_2,e_3}$, it defines an element 
of the 2-torsion part of $WC(E_{e_1,e_2,e_3},\mathbb{Q})$ (cf. [20], Ex. 4.5.1). \smallskip 

The assignment $(x_0,x_1,x_2,x_3) \mapsto \bigl( (b_1x_1^2/x_0^2)+e_1, \, b_1b_2x_1x_2x_3/x_0^3
\bigr)$ extends to a well-defined regular mapping $Q_{e_1,e_2,e_3,b_1,b_2} \rightarrow 
E_{e_1,e_2,e_3}$ of degree 4. In particular, if $(x_0,x_1,x_2,x_3)$ is a rational point on 
the homogeneous space $Q_{e_1,e_2,e_3,b_1,b_2}$, one gets a rational point on the given curve. 
Note that the logarithmic height of the point $(b_1x_1^2/x_0^2+e_1,b_1b_2x_1x_2x_3/x_0^3)$ 
on $E_{e_1e_2e_3}$ is about three times the height of the point $(x_0,x_1,x_2,x_3)$ on 
$Q_{e_1,e_2,e_3,b_1,b_2}$. (Recall that the logarithmic height of a point 
$P\in{\mathbb{P}}^N(\mathbb{Q})$ is the logarithm of $\max (|x_0|,\ldots ,|x_N|)$ where 
$P=(x_0:x_1:\cdots :x_N)$ is represented with coprime integers $x_i$.) This explains why it 
is reasonable to look for 
points on $Q_{e_1,e_2,e_3,b_1,b_2}$ instead of finding rational points on $E_{e_1,e_2,e_3}$ 
directly. The problem is to determine the approriate homogeneous spaces which have 
rational points. The most important approach to this question is contained in the following 
two-descent procedure (cf. [20], Chap. X, Prop. 1.4; [8], Chap. IV, Section 3). \\ 

Let $S$ be the set of integers consisting of $-1$, $2$ and all the divisors of the discriminant 
of $E_{e_1,e_2,e_3}$. Let $\mathbb{Q}(S,2)$ be the subgroup of $\mathbb{Q}^{*}/(\mathbb{Q}^{*}
)^{2}$ generated by the elements of $S$. Thus any element in $\mathbb{Q}(S,2)$
has a unique representation by an integer which is squarefree and which has prime factors 
only in $S$. Let us consider the mappings $\varphi_{i}:E_{e_1,e_2,e_3}(\mathbb{Q})/2E_{
e_1,e_2,e_3}(\mathbb{Q})\rightarrow\mathbb{Q}(S,2)$ defined by
$$\varphi_{i}(P):=\begin{cases}
   x-e_{i} & \hbox{if}\ P=(x,y)\ \hbox{with}\ x\neq e_i; \\
   (e_i-e_j)(e_i-e_k) & \hbox{if}\ P=(e_i,0)\ \hbox{where}\ \{ i,j,k\} = \{ 1,2,3\}; \\
   1 & \hbox{if}\ P=\infty. 
\end{cases}\leqno{(10)}$$
Here the values are taken as representatives modulo $(\mathbb{Q}^{*})^{2}$.
Then $\varphi_{i}$ is well-defined and a homomorphism of groups. Furthermore, 
the homomorphism
$$\varphi: 
  \begin{matrix} 
  E_{e_1,e_2,e_3}(\mathbb{Q})/2E_{e_1,e_2,e_3}(\mathbb{Q}) & \rightarrow & \mathbb{Q}(S,2)
      \times\mathbb{Q}(S,2) \\ 
 P & \mapsto & \bigl(\varphi_1(P), \varphi_2(P)\bigr)
\end{matrix}\leqno{(11)}$$
is injective (see [8], Prop. 4.8, or [20], Chap. X, Prop. 1.4).\\
\\
{\bf Consequence:} The homogeneous spaces $Q_{e_1,e_2,e_3,b_1,b_2}$ which contain a rational 
point (and hence are trivial in the sense of [20], Chap. X, Section 3) are exactly those for 
which $(b_1,b_2)$ is contained in the image of the mapping $\varphi$. Since $\mathbb{Q}(S,2)
\times\mathbb{Q}(S,2)$ is finite, we therefore have to check only a finite number of homogeneous 
spaces for the existence of a rational point.\\
\\
{\bf Remark:} If the curve is given in the Weierstra\ss{} form for $E_{e_1,e_2,e_3}$, the two-descent 
procedure (and thus the meaning of the mapping $\varphi$) can be explained in a very elementary 
way. The affine part of $E_{e_1,e_2,e_3}$ is given by the equation $y^2=(x-e_1)(x-e_2)(x-e_3)$. 
If $p=(x,y)\in\mathbb{Q}^2$ is a point on $E_{e_1,e_2,e_3}(\mathbb{Q})$ such that all three 
factors $(x-e_i)$ are rational squares, then $p=[2]\cdot q$ for some rational point $q\in E_{
e_1,e_2,e_3}(\mathbb{Q})$, where $[2]$ denotes the multiplication-by-two map on $E_{e_1,e_2,e_3}$.
This halving procedure stops after a finite number of steps (because $E_{e_1,e_2,e_3}(\mathbb{Q})$ 
is finitely generated), yielding a point $p\in E_{e_1,e_2,e_3}(\mathbb{Q})$ which is not twice 
another rational point on $E_{e_1,e_2,e_3}(\mathbb{Q})$. For this point $p=(x,y)$ we know that 
not all the factors $(x-e_i)$ can be rational squares, whereas the product $y^2=(x-e_1)(x-e_2)
(x-e_3)$ is a rational square. Writing 
$$x-e_1=A_1\alpha_1^2,\qquad x-e_2=A_2\alpha_2^2,\qquad x-e_3=A_3\alpha_3^2\leqno{(12)}$$ 
with squarefree integers $A_i$, we can eliminate $x=A_1\alpha_1^2+e_1$ to arrive at 
the equations
$$A_1\alpha_1^2-A_2\alpha_2^2 = e_2-e_1,\qquad A_1\alpha_1^2-A_3\alpha_3^2=e_3-e_1.\leqno{(13)}$$
Moreover, the product $A_1A_2A_3$ is a perfect square, which implies that $A_3=A_1A_2$ 
up to a square factor. In this way we obtain the equations of the homogeneous space in 
the above considerations. \medskip 

\subsubsection{Reducing the number of possibilities}

Let us first observe that, given a set $P\subseteq E_{e_1,e_2,e_3}(\mathbb{Q})$ of known 
rational points on $E_{e_1,e_2,e_3}$, and given a homogeneous space $Q_{e_1,e_2,e_3,b_1,b_2}$ 
by coefficients $(b_1,b_2)\in\mathbb{Q}^{*}/(\mathbb{Q}^{*})^{2}\times\mathbb{Q}^{*}/(
\mathbb{Q}^{*})^{2}$, then $Q_{e_1,e_2,e_3,b_1,b_2}$ contains a rational point (i.e., is 
trivial in the Weil-Chatelet group) if and only if for any point $p\in P$ the homogeneous 
space $Q_{e_1,e_2,e_3,c_1,c_2}$ with $(c_1,c_2)=(b_1,b_2)\cdot\varphi(p)$ has a rational point.
So we may say that $(b_1,b_2),(c_1,c_2)\in\mathbb{Q}^{*}/(\mathbb{Q}^{*})^{2}\times
\mathbb{Q}^{*}/(\mathbb{Q}^{*})^{2}$ are $P$-equivalent if and only if there is a point 
$p\in P$ such that $(c_1,c_2)=(b_1,b_2)\cdot\varphi(p)$. Using this relation, we may 
restrict our considerations to special representatives of homogeneous spaces whenever we 
{\it a priori} know some rational points on $E_{e_1,e_2,e_3}$. Note that a well-known set 
of rational points on $E_{e_1,e_2,e_3}$ is always given by the set of 2-torsion points.\\
\\
Now we consider either a single curve $E_{e_1,e_2,e_3}$ or else a family of curves 
$E_{e_1,e_2,e_3}$ depending on some parameter(s) such that we know in advance that the 
curve(s) contain rational points other than the 2-torsion points. To find appropriate
homogeneous spaces for $E_{e_1,e_2,e_3}$ possessing rational points, we may pursue the 
following strategy:
\begin{itemize}
\item determine the (finite) set $\mathbb{Q}(S,2)$; 
\item build equivalence classes of pairs $(b_1,b_2)\in\mathbb{Q}(S,2)\times\mathbb{Q}(S,2)$ 
  which are $P$-equivalent with respect to the set $P$ of 2-torsion points of $E_{e_1,e_2,e_3}$; 
\item exclude those classes of pairs $(b_{1},b_{2})\in\mathbb{Q}(S,2)\times\mathbb{Q}(S,2)$
  for which a rational solution cannot exist for one of the following reasons:
  \begin{itemize}
  \item the quadric intersection $Q_{e_1,e_2,e_3,b_1,b_2}$ has no rational points because of  
    obstructions due to properties of quadratic residues; 
\item at least one of the quadratic equations defining $Q_{e_1,e_2,e_3,b_1,b_2}$ has no 
  solution (for example because of Legendre's criterion). 
\end{itemize}
\end{itemize}
~\\
Then the homogeneous spaces associated with the remaining classes are good candidates for 
having rational points and thus giving rational points on the original elliptic curve(s) 
$E_{e_1,e_2,e_3}$. Examples for the application of this method will be presented in the 
next section.

\subsubsection{Notations }

The notations used up to this point were chosen to be consistent with the ones used in 
[20] and [8]. For the use in later sections, it will be convenient to modify the notations 
slightly. As mentioned in 2.1, we always may assume that the elliptic curve under 
consideration is one of the curves $E_{M,N}$ with different nonzero integers $M,N\in\mathbb{Z}
\setminus\{ 0\}$, given in affine form by a Weierstra\ss{} equation of the special form $y^2 = 
x(x+M)(x+N)$. In addition, we may assume $M=pk>0$ and $N=-qk<0$ with coprime natural numbers 
$p,q\in\mathbb{N}$ and a square-free natural number $k$. We write the equations for the 
homogeneous spaces in the form 
$$AX_0^2+MX_1^2-BX_2^2=0,\qquad AX_0^2+NX_1^2-CX_3^2=0\leqno{(14)}$$ 
where $A,B,C\in\mathbb{Q}(S,2)$ are represented by squarefree integers, with $C$ depending on 
$A$ and $B$ by the condition that $ABC$ is a perfect square. We call $(A,B,C)$ a triplet defining 
a homogeneous space for $E_{M,N}$. The approach sketched in 2.3.2 provides us with means to 
reduce the number of potential triplets $(A,B,C)$ by ruling out impossible ones and by 
identifying triplets with respect to 2-torsion equivalence. (In the next section some speficic 
examples are provided which explicitly show how this is done.) \smallskip 

Let us choose any of the possible 2-descent parameter sets 
$(A,B,C)$. Then we know that finding the expected solution of $E_{M,N}$ is 
equivalent to finding a solution of the system $(12)$, which in our situation 
reads
$$x=A\alpha^2, \qquad x+M=B\beta^2, \qquad x+N=C\gamma^2.\leqno{(15)}$$
where $x,\alpha,\beta,\gamma\in\mathbb{Q}$. With such a solution
we get a rational point $(x,y)$ on $E_{M,N}$ by setting $y=\sqrt{ABC}\alpha\beta\gamma$. 
The system (15) is equivalent to the system (14), and an integer solution $(x_0,x_1,x_2,x_3)$ 
of (14) yields a rational solution of (15) via 
$$x = \frac{Ax_0^2}{x_1^{2}}, \qquad y=\sqrt{ABC}\,\frac{x_0x_2x_3}{x_1^3}, \qquad 
  (\alpha,\beta,\gamma)=\left(\frac{x_0}{x_1},\frac{x_2}{x_1},\frac{x_3}{x_1}\right). 
  \leqno{(16)}$$ 
In projective notation, the sought rational point $(x,y)$ on $E_{M,N}$ is given by the 
expression $(x_1^3 : Ax_0^2x_1 : \sqrt{ABC}\,x_0x_2x_3)$. We see that the logarithmic height 
of this point is about 3 times the logarithmic height of the solution $(x_0,x_1,x_2,x_3)$.
\smallskip 

In the sequel we thus may restrict our considerations to systems of quadratic equations 
of the form (14). We may also assume that this system has a solution at all, which in 
particular implies that each of the two individual quadratic equations has a solution.
 
\section{Examples for the two-descent}

As before we consider elliptic curves which are given by an affine equation 
of the form $y^2=x(x+M)(x+N)$ where $M=pk$ and $N=-qk$ with coprime natural numbers $p,q$ and 
a squarefree natural number $k$. We denote the elliptic curve defined 
in this way by $E_{p,q,k}$ and its Mordell-Weil group by $E_{p,q,k}(\mathbb{Q})$.
As explained before, our aim is to determine the homogeneous spaces
over $E_{p,q,k}$ defined over the rationals which have a rational point. 
For this purpose we consider the homomorphism $\varphi=(\varphi_{-pk},\varphi_{0},
\varphi_{qk}):E_{p,q,k}(\mathbb{Q})/2E_{p,q,k}(\mathbb{Q})\rightarrow(\mathbb{Q}^{*}/
(\mathbb{Q}^{*})^{2})^{3}$ where the components are defined as in (11). For the
sake of simplicity we look at all three components even though the third
one is already determined by the other two. From the considerations in 2.3.1 
above we know that finding non-2-torsion points on $E_{p,q,k}$ is
equivalent to finding triplets in the image of the mapping $\varphi$ 
(or, equivalently, homogeneous spaces having a rational point). In this
section we give examples how to find candidates of triplets which potentially 
may lie in the image of $\varphi$. In other words, we try to find criteria which exclude 
many of these triplets from potentially lying in the image of $\varphi$. 
\medskip 

\subsection{General strategy}

We proceed as follows:
\begin{enumerate}
\item Determine the discriminant of the equation for $E_{p,q,k}$. Determine
  all prime divisors of this discriminant and the set $S_{p,q,k}(2,\mathbb{Q})$
  of all squarefree integers which are candidates for components appearing
  in the image of $\varphi$.
\item Determine all the triplets $(A,B,C)\in(S_{p,q,k}(2,\mathbb{Q}))^3$ 
  which cannot {\it a priori} be ruled out to lie in the image of $\varphi$. Note that
  with the conventions set up above, $A$ is positive and $ABC$ is a perfect square.
\item Generate equivalence classes of the remaining triplets with respect
  to 2-torsion equivalence; extract representatives of these equivalence classes.
\item Eliminate all those triplets which give rise to a homogeneous space over $E_{p,q,k}$ 
  which cannot have a rational point due to the quadratic residue behaviour of the quadratic
  equations defining the homogeneous space.
\item Eliminate all triplets leading to a system of quadratic equations
  of the form (14) of section 2 such that at least one of the quadratic
  equations is not solvable.
\end{enumerate}
The remaining triplets are good candidates to start the algorithm to be developed in the sequel.\\

\subsection{General data}

Since the discriminant of the elliptic curve $E_{p,q,k}$ is given by $\hbox{discr}(E_{p,q,k}) = 
16p^2q^2(p+q)^2k^6$, the elements of $S_{p,q,k}(2,\mathbb{Q})$ are the squarefree
integral numbers which are composed of $-1$, $2$ and the prime divisors of $p,q,p+q$ and 
$k$. For analyzing the 2-torsion-equivalent triplets, it is interesting
to know the values of the mappings $\varphi_{-pk}$, $\varphi_{0}$,
$\varphi_{qk}$ at the 2-torsion points $(-pk,0)$, $(0,0)$, $(qk,0)$.
We collect the results in the following table, where the values always
have to be regarded as representatives which should be replaced by
their squarefree parts.
$$\begin{tabular}{|c||c|c|c|}
  \hline 
   & $(-pk,0)$ & $(0,0)$ & $(qk,0)$\tabularnewline
  \hline 
  \hline 
  $\varphi_{-pk}$ & $p(p+q)$ & $pk$ & $(p+q)k$\tabularnewline  
\hline 
$\varphi_{0}$ & $-pk$ & $-pq$ & $qk$\tabularnewline
\hline 
$\varphi_{qk}$ & $-(p+q)k$ & $-qk$ & $q(p+q)$\tabularnewline
\hline 
\end{tabular}$$ 
For any triplet $(A,B,C)$ defining a homogeneous space over $E_{p,q,k}$
we consider the equations $x+pk=A\alpha^2$, $x=B\beta^2$, $x-qk=C\gamma^2$  
and try to examine whether or not these equations are solvable in rational
numbers $x,\alpha,\beta,\gamma$. Here we may restrict our attention
to some representative $(A,B,C)$ modulo 2-torsion equivalence.

\subsection{Examples}

\subsubsection{Congruent prime numbers}

The case of congruent prime numbers corresponds to the case $p=q=1$ with 
$k$ being a prime number. This case is discussed extensively in [8]. The 
results can be summarized as follows:
\begin{itemize}
\item If $k\equiv1$ or $3$ mod $8$, then the rank of $E_{1,1,k}(\mathbb{Q})$
  is zero and there is no non-2-torsion point on $E_{1,1,k}(\mathbb{Q})$
\item If $k\equiv5$ mod $8$, then the triplets $(A,B,C)$ yielding a homogeneous
  space with a rational point are given by the 2-torsion-equivalence
  class containing $(1,-1,-1)$, $(2,k,2k)$, $(k,1,k)$ and $(2k,-k,-2)$. 
\item If $k\equiv 7$ mod $8$, then the triplets $(A,B,C)$ yielding a homogeneous
  space with a rational point are given by the 2-torsion-equivalence
  class containing $(2,1,2)$, $(1,-k,-k)$, $(2k,-1,-2k)$ and $(k,k,1)$. 
\end{itemize}

\subsubsection{Congruent numbers which are twice a prime number}

This case corresponds to $p=q=1$ and $k=2\ell$ where $\ell$ is a prime number and can 
be discussed in a way similar to the above case. By different methods 
(see [18], p. 343) we {\it a priori} know that $2\ell$ is a congruent number 
if $\ell\equiv 3$ or $7$ mod $8$. In this situation we have $S_{1,1,2\ell} = \{-1,2,\ell\}$, 
and the possibilities for the triplets defining homogeneous spaces over $E_{1,1,2\ell}$ are 
grouped into the following 2-torsion-equivalence classes:\smallskip 

1. $\{(1,1,1),(2,-2l,-l),(2l,-1,-2l),(l,2l,2)\}$ \hfill\break\indent 
2. $\{(1,-1,-1),(2,2l,l),(2l,1,2l),(l,-2l,2)\}$ \hfill\break\indent 
3. $\{(1,2,2),(2,-l,-2l),(2l,-2,-l),(l,l,1)\}$ \hfill\break\indent 
4. $\{(1,-2,-2),(2,l,2l),(2l,2,l),(l,-l,-1)\}$ \hfill\break\indent 
5. $\{(2,1,2),(1,-2l,-2l),(l,-1,-l),(2l,2l,1)\}$ \hfill\break\indent 
6. $\{(2,-1,-2),(1,2l,2l),(l,1,l),(2l,-2l,-1)\}$ \hfill\break\indent 
7. $\{(2,2,1),(1,-l,-l),(l,-2,-2l),(2l,l,2)\}$ \hfill\break\indent 
8. $\{(2,-2,-1),(1,l,l),(l,2,2l),(2l,-l,-2)\}$ 
\smallskip 

\noindent 
The first of these classes is given by the 2-torsion elements themselves.
For the other ones it is sufficient to consider the equations
corresponding to the first triplet.\\
\\
2. If $(A,B,C)=(1,-1,-1)$ then $x+2\ell = \alpha^2$, $x = -\beta^2$ and $x-2\ell = -\gamma^2$.  
Subtracting the third equation from the first yields $4\ell = \alpha^2 + \gamma^2$. 
Clearing denominators and reducing modulo the prime number $\ell$, we see that 
$-1$ is a quadratic residue mod $\ell$, hence $\ell\equiv 1$ mod $4$. So $\ell\equiv 1$ 
or $5$ mod $8$.\\
\\
3. If $(A,B,C)=(1,2,2)$ then $x+2\ell = \alpha^2$, $x=2\beta^2$ and $x-2\ell = 2\gamma^2$.  
Subtracting the third equation from the first yields $4\ell = \alpha^2-2\gamma^2$. 
Clearing denominators and reducing modulo $\ell$ we see that $2$ is a quadratic residue 
mod $\ell$, hence $\ell\equiv 1$ or $7$ mod $8$.\\
\\
4. If $(A,B,C)=(1,-2,-2)$ then $x+2\ell = \alpha^2$, $x=-2\beta^2$ and $x-2\ell = -2\gamma^2$.  
Subtracting the second equation from the first yields $2l\ell = \alpha^2 + 2\beta^2$. 
Clearing denominators and reducing modulo $\ell$ we see that $-2$ is a quadratic residue mod $\ell$, 
hence $\ell\equiv 1$ or $3$ mod $8$.\\
\\
5. If $(A,B,C)=(2,1,2)$ then $x+2\ell = 2\alpha^2$, $x=\beta^2$, $x-2\ell = 2\gamma^2$.  
Subtracting the second equation from the first yields $2\ell = 2\alpha^2-\gamma^2$. 
Clearing denominators and reducing modulo $\ell$ we see that $-2$ is a quadratic residue mod 
$\ell$, hence $\ell\equiv 1$ or $7$ mod $8$.\\
\\
6. If $(A,B,C)=(2,-1,-2)$ then $x+2\ell = 2\alpha^2$, $x=-\beta^2$, $x-2\ell = -2\gamma^2$.  
On the one hand, subtracting the third equation from the first yields 
$4\ell = 2\alpha^2 + 2\gamma^2$. Clearing denominators, dividing by $2$ and reducing modulo $\ell$ 
we see that $-1$ is a quadratic residue mod $\ell$, hence $\ell\equiv 1$ mod $4$ so that $\ell\equiv 
1$ or $5$ mod $8$. On the other hand, subtracting the second equation from the first 
yields $2\ell = 2\alpha^2 + \beta^2$. Clearing denominators and reducing modulo $\ell$ we see that 
$-2$ is a quadratic residue mod $\ell$, which implies that $\ell\equiv 1$ or $3$ mod $8$. 
Combining both results, we see that $\ell\equiv 1$ mod $8$.\\
\\
7. If $(A,B,C)=(2,2,1)$ then $x+2\ell = 2\alpha^2$, $x = 2\beta^2$, $x-2\ell = \gamma^2$.  
Subtracting the third equation from the first yields $4\ell = 2\alpha^2 - \gamma^2$. 
Clearing denominators and reducing modulo $\ell$ we see that $2$ is a quadratic residue 
mod $\ell$, hence $\ell\equiv 1$ or $7$ mod $8$.\\
\\
8. If $(A,B,C)=(2,-2,-1)$ then $x+2\ell = 2\alpha^2$, $x = -2\beta^2$, $x-2\ell = -\gamma^2$.  
On the one hand, subtracting the third equation from the first yields $4\ell = 
2\alpha^2 + \gamma^2$. Clearing denominators and reducing modulo $\ell$ we see 
that $-2$ is a quadratic residue mod $\ell$, hence $\ell\equiv 1$ or $3$ mod $8$. 
On the other hand, subtracting the second equation from the first yields 
$2\ell = 2\alpha^2 + 2\beta^2$. Clearing denominators, dividing by 2 and reducing 
modulo $\ell$ we see that $-1$ is a quadratic residue mod $\ell$, hence $\ell\equiv 1$ 
or $5$ mod $8$. Combining both results we see that $\ell\equiv 1$ mod $8$. \\
\\
Since we know beforehand (see [18], p. 343) that the number $2\ell$ is congruent 
if $\ell\equiv 3$ or $7$ mod $8$, the above calculations give the following
necessary conditions. 
\begin{itemize}
\item If $l\equiv3$ mod $8$ then the only possibilities for solvable systems
are given by the equivalence classes of the triplets $(1,-2,-2)$
and $(2,-2,-1)$.
\item If $l\equiv7$ mod $8$ then the only possibilities for solvable systems
are given by the equivalence classes of the triplets $(1,2,2)$, $(2,1,2)$
and $(2,2,1)$.
\end{itemize}
Explicit calculations show that in the first case the triplet $(1,-2,-2)$ yields 
solvable equations whereas the second one does not. Also by explicit calculations, 
in the second case the only triplet yielding solvable equations is $(2,2,1)$.
For the rational points on the corresponding elliptic curves we refer
to section 5 below. 

\subsubsection{Examples with rank $\bigl(E_{p,q,k}(\mathbb{Q})\bigr)=2$}

The examples considered here correspond to the $2\pi/3$-congruent
numbers $14$, $206$ and $398$. These numbers are of the form $2\ell$ 
where $\ell$ is a prime number satisfying $\ell\equiv 7$ mod $96$; hence 
the associated elliptic curve is given by $p=1$, $q=3$ and $k=2\ell$. 
We do not know whether all of these numbers $\ell$ yield elliptic curves 
with $\hbox{rank}(E_{1,3,2\ell}(\mathbb{Q}))=2$, but at least for the 
three explicit examples above we will show in section 5 that there are 
two independent solutions.\\

The discriminant of $E_{1,3,2\ell}$ is $\hbox{\rm discr}(E_{1,3,2\ell}) = 2^{14}3^2\ell^{6}$ 
so that $S_{1,3,2\ell}=\{-1,2,3,\ell\}$. In particular, we get $128$
possibilities for the triplets $(A,B,C)$ characterizing the homogeneous
spaces associated with the 2-descent. The table for the mapping $\varphi$
at the 2-torsion points looks as follows. 
$$\begin{tabular}{|c||c|c|c|}
  \hline 
   & $(-2l,0)$ & $(0,0)$ & $(6l,0)$\tabularnewline
  \hline 
  \hline 
  $\varphi_{-2l}$ & $1$ & $2l$ & $2l$\tabularnewline
  \hline 
  $\varphi_{0}$ & $-2l$ & $-3$ & $6l$\tabularnewline
  \hline 
  $\varphi_{6l}$ & $-2l$ & $-6l$ & $3$\tabularnewline
  \hline 
  \end{tabular}$$ 
We see that any 2-torsion equivalence class contains a triplet 
which is independent of $2\ell$. So the equivalence classes of these
triplets are already determined by the 32 possible triplets composed of
the factor set $\{-1,2,3\}$. (Note that $-1$ does not occur as a factor 
of $A$ since $A$ is positive.) \smallskip 

Let us discuss the explicit example of the triplet $(A,B,C)=(3,3,1)$.
The equations from 3.2 are $x+2\ell = 3\alpha^2$, $x=3\beta^2$ and 
$x-6\ell = \gamma^2$. Subtracting the third equation from the first yields 
$8\ell = 3\alpha^2 - \gamma^2$. Clearing denominators and reducing modulo 
the prime number $\ell$ shows that $3$ is a quadratic residue mod $\ell$, 
hence $\ell\equiv 1$ or $11$ mod $12$. But this is not the case for $l\equiv 7$ 
mod $96$. So this triplet does not belong to a homogeneous space having a rational
point. \smallskip 

In the same way we may discuss all 32 possibilities. The only triplets which do not 
lead to a contradiction are $(1,2,2)$, $(1,-3,-3)$, $(1,-6,-6)$, $(2,1,2)$, $(2,2,1)$, 
$(2,-3,-6)$ and $(2,-6,-3)$. These triplets provide candidates for the search for 
rational solutions of the corresponding equations. In section 5 we will see that the 
triplets $(1,2,2)$, $(2,-3,-6)$ and $(2,-6,-3)$ indeed lead to solutions, as, of course, 
also do all triplets which are 2-torsion-equivalent to one of these. For the sake of 
completeness, let us explicitly write down the 2-torsion equivalence classes of these 
three triplets: \smallskip 
 
$\bullet$ $\{(1,2,2),(2l,-6,-3l),(1,-l,-l),(2l,3l,6)\}$; \hfill\break\indent  
$\bullet$ $\{(2,-3,-6),(l,1,l),(2,6l,3l),(l,-2l,-2)\}$; \hfill\break\indent 
$\bullet$ $\{(2,-6,-3),(l,2,2l),(2,3l,6l),(l,-l,-1)\}$. \smallskip 

Note that in the group structure of $(\mathbb{Q}^{*}/(\mathbb{Q}^{*})^{2}\times(\mathbb{Q}^{*}/
(\mathbb{Q}^{*})^{2}\times(\mathbb{Q}^{*}/(\mathbb{Q}^{*})^{2}$ the third triplet is the product 
of the other two. So solutions corresponding to the third triplet correspond to sums of solutions
of the other ones (the sum being calculated in terms of the elliptic
curve structure of the given curve $E_{1,3,2\ell}$). But the solutions
belonging to the first and the second triplet are independent with
respect to the elliptic curve addition, and since they are of infinite
order we have shown that the Mordell-Weil rank of these groups is
at least $2$. 
 
\section{Algorithm}

\subsection{Weak form of the algorithm }

The starting point of our algorithm is a system of quadratic equations
of the form\\
\parbox{.05\textwidth}{(17)} 
\parbox{.95\textwidth}{ 
\begin{alignat*}{3} 
  & Q_1(X_0,X_1,X_2)\ &&=\ a_{00}X_0^2 + a_{11}X_1^2 + a_{22}X_2^2\ &&=\ 0 \\
  & Q_2(X_0,X_1,X_3)\ &&=\ b_{00}X_0^2 + b_{11}X_1^2 + b_{33}X_3^2\ &&=\ 0
\end{alignat*}} \\
which we assume to be solvable in integers. Then in particular there
is a point $(x_0,x_1,x_2)\in\mathbb{Z}^{3}\setminus\{(0,0,0)\}$ with 
$Q_1(x_0,x_1,x_2)=0$. By Newton's method, projecting from $(x_0,x_1,x_2)$,  
we may parametrize the solutions of $Q_1$ by quadratic polynomials
$$X_0=\Phi_{0}(\xi_0,\xi_1), \qquad X_1=\Phi_1(\xi_0,\xi_1), \qquad 
  X_2=\Phi_{2}(\xi_0,\xi_1), \leqno{(18)}$$ 
which means that the solutions $(x_0,x_1,x_2)$ of the equation $Q_1=0$ are exactly 
the points $\bigl(\Phi_0(\xi_0,\xi_1), \Phi_1(\xi_0,\xi_1), \Phi_2(\xi_0,\xi_1)\bigr))$
where $(\xi_0,\xi_1)\in\mathbb{P}^{1}(\mathbb{Q})$. Substituting this parametrization 
for $X_0$ and $X_1$ into the equation for $Q_2$ yields an equation\\
\parbox{.05\textwidth}{(19)} 
\parbox{.95\textwidth}{
  \begin{alignat*}{2}
  & && Q_3(\xi_0,\xi_1,X_3)\ =\ Q_2\bigl(\Phi_0(\xi_0,\xi_1), \Phi_0(\xi_0,\xi_1), X_3\bigr)\\ 
  & && =\ b_{00}\Phi_{0}(\xi_0,\xi_1)^2 + b_{11}\Phi_1(\xi_0,\xi_1)^2 + b_{33}X_3^2\ =\ 0
  \end{alignat*} }\\  
which is of degree $4$ in $(\xi_0,\xi_1)$ and of degree $2$ in $X_3$. This equation is the 
basis of the weak form of the following algorithm for finding solutions of a system of the 
form (17). \\
\\
\def\eins{\phantom{eee}} \def\zwei{\eins\eins} \def\drei{\zwei\eins} \def\vier{\drei\eins}
\eins\texttt{/{*} weak algorithm {*}/} \hfill\break
\eins\texttt{INPUT: quadrics $Q_1$ and $Q_2$} \hfill\break
\eins\texttt{BEGIN} \hfill\break
\eins\texttt{determine a solution $(x_0,x_1,x_2)$ of $Q_1=0$} \hfill\break
\eins\texttt{calculate parametrizations $X_i=\Phi_i(\xi_0,\xi_1)$, $i=0,1,2$} \hfill\break
\eins\texttt{loop over $(\xi_0,\xi_1)$} \hfill\break
\zwei\texttt{BEGIN} \hfill\break
\drei\texttt{calculate value$=-b_{33}(b_{00}\Phi_0(\xi_0,\xi_1)^2+b_{11}\Phi_1(\xi_0,\xi_1)^2)$} 
     \hfill\break
\drei\texttt{check whether value is a square number} \hfill\break
\drei\texttt{if yes} \hfill\break
\drei\texttt{BEGIN} \hfill\break
\vier\texttt{/{*} solution found {*}/} \hfill\break
\vier\texttt{$x_0=\Phi_0(\xi_0,\xi_1)$} \hfill\break
\vier\texttt{$x_1=\Phi_1(\xi_0,\xi_1)$} \hfill\break
\vier\texttt{$x_2=\Phi_2(\xi_0,\xi_1)$} \hfill\break 
\vier\texttt{$x_3=\sqrt{\hbox{\texttt value}}/(-b_{33})$} \hfill\break
\vier\texttt{RETURN $(x_0,x_1,x_2,x_3)$} \hfill\break
\drei\texttt{END} \hfill\break
\zwei\texttt{END} \hfill\break
\eins\texttt{END} \hfill\break
\\
{\bf Remarks:} 
\begin{itemize} 
\item The algorithm is not guaranteed to terminate. In fact, to the best of our knowledge 
   no upper bound for the smallest solution of a system of the form (17) is known.
\item The logarithmic height of a solution $(x_0,x_1,x_2,x_3)$ found by the algorithm is about 
  twice the logarithmic height of the parameters $(\xi_0,\xi_1)$. 
\item Assuming that the algorithmic complexity of the calculations within the algorithm do not 
  too strongly depend on the magnitude of the numbers involved, then the complexity of the 
  algorithm depends quadratically on the height of the parameters $(\xi_0,\xi_1)$. 
\item Using a general purpose computer, it is possible to calculate the 
  final loop for the parameters $\xi_0,\xi_1$ up to about $10\, 000$ within a few minutes, which  
  yields results for the original quadratic equations with up to about $25$ digits (in decimal 
  representation).
\end{itemize} 

\subsection{Strong form of the algorithm}

This section is devoted to an improved version of the weak form of the
algorithm, provided a special property of the considered equations holds. 

\subsubsection{Additional preparation and condition}

We note that the starting point of our algorithm is the system $(17)$ which 
consists of two quadratic equations of a very special form. Namely, exactly two of the 
variables appear in both equations. We can, of course, eliminate from each of the equations 
one of the common variables, which yields four equations each of which contains exactly 
three of the four original variables, and each choice of two of these equations defines the 
same quadric intersection in projective 3-space. We now impose a condition which seems to be 
of a purely technical nature (not corresponding to any geometric or arithmetic property of the
given elliptic curve), but seems to be satisfied in many particular cases. (See the examples 
in the subsequent sections.) \\
\\
{\bf Condition:} We suppose that at least one of the four equations has
a solution with one coordinate being zero. \\
\\
We choose the equation of the condition to be the first one and choose
a second equation such that the variable whose coordinate in the distinguished
solution is zero appears in both equations. After renaming the variables if necessary, 
we may suppose without loss of generality that the equations have the form $(17)$ and 
that $(0,x_1,x_2)$ is an integer solution for $Q_1$. As above, we parametrize the quadratic 
form $Q_1$ projecting from this point $(0,x_1,x_2)$, which gives expressions
$X_i=\Phi_i(\xi_0,\xi_1)$ as in $(18)$, and as above we substitute the parametrizations 
of $X_0$ and $X_1$ into the second equation. The crucial observation is that in our 
particular situation the polynomial $Q_3$ in equation $(19)$ is biquadratic in $(\xi_0,\xi_1)$; 
see 2.2.3. 
Hence substituting $Y_0=\xi_0^2$, $Y_1=\xi_1^2$, $Y_2=X_3$ yields a quadratic form, which we 
denote (with a slight abuse of notation) by $Q_3(Y_0,Y_1,Y_2)$. Our task is then to find an 
integer solution of the equation $Q_3(y_0,y_1,y_2)=0$ in such a way that both $y_0$ and $y_1$ 
are squares. \\
\\
{\bf Remark:} The property of being a square number is not invariant with
respect to the change of representatives of projective coordinates.
So the above problem can be stated equivalently by saying that we
look for a coprime solution $(y_0,y_1,y_2)$ for $Q_3$ and a factor $\mu$ 
(which we may assume to be squarefree) such that $y_0=\mu\sigma_0^2$,
$y_1=\mu\sigma_1^2$. This seems to be difficult. However, we can establish a strong 
condition for the possible factors $\mu$, as will be explained below.\\
\\
{\bf Remark:} The order of magnitude of $\xi_{i}$ is about twice the order
of magnitude of $X_{i}$, so $Y_{i}$ and $X_{i}$ are nearly of the
same order of magnitude.

\subsubsection{Parametrization of $Q_{3}$ and condition for $\mu$}

As always we assume our problem to have a solution, hence $Q_3$ has a solution with the 
required property; in particular, $Q_3$ has a rational solution which, without loss of 
generality, we may assume to consist of coprime integers. Let $p=(p_0,p_1,p_2)$ be such 
a solution. We parametrize $Q_3$ using $p$ as projection point according to Newton's method to get
representations
$$Y_0\ =\ \Psi_0(\eta_0,\eta_1), \qquad 
  Y_1\ =\ \Psi_1(\eta_0,\eta_1), \qquad 
  Y_2\ =\ \Psi_2(\eta_0,\eta_1) \leqno{(20)}$$ 
with quadratic forms $\Psi_0,\Psi_1,\Psi_2$. We observe that if $\mu$ is a factor as considered 
in the last remark of 4.2.1, then $\mu$ divides the numbers $\Psi_i(\eta_0,\eta_1)$ where 
$i=0,1,2$ and hence also divides all ${\mathbb Z}$-linear combination of these three numbers. 
We write\\
\parbox{.05\textwidth}{(21)} 
\parbox{.95\textwidth}{
  \begin{eqnarray*}
  \Psi_0(\eta_0,\eta_1) &=& \psi_{00}^{(0)}\eta_0^2 + \psi_{01}^{(0)}\eta_0\eta_1 
        + \psi_{11}^{(0)}\eta_1^2 \\
  \Psi_1(\eta_0,\eta_1) &=& \psi_{00}^{(1)}\eta_0^2 + \psi_{01}^{(1)}\eta_0\eta_1 
        + \psi_{11}^{(1)}\eta_1^2 \\
  \Psi_2(\eta_0,\eta_1) &=& \psi_{00}^{(2)}\eta_0^2 + \psi_{01}^{(2)}\eta_0\eta_1 
        + \psi_{11}^{(2)}\eta_1^2 
  \end{eqnarray*}} \\  
and determine coprime integer coefficients $c_i$ where $i=0,1,2$ such that\\ 
\parbox{.05\textwidth}{(22)} 
\parbox{.95\textwidth}{
  \begin{eqnarray*} 
  c_0\psi_{00}^{(0)} + c_1\psi_{00}^{(1)} + c_2\psi_{00}^{(2)} &=& 0, \\ 
  c_0\psi_{11}^{(0)} + c_1\psi_{11}^{(1)} + c_2\psi_{11}^{(2)} &=& 0.
  \end{eqnarray*}} \\  
so that  
$$c_0Y_0 + c_1Y_1 + c_2Y_2\ =\ (c_0\psi_{01}^{(0)} + c_1\psi_{01}^{(1)} + c_2
  \psi_{01}^{(2)})\cdot\eta_0\eta_1. \leqno{(23)}$$ 
Since we are looking for coprime solutions $(y_{0},y_{1},y_{2})$
we may assume $(\eta_{0},\eta_{1})$ to be coprime and the factor
$\mu$ to be coprime with both $\eta_0$ and $\eta_1$. So we have
the following necessary condition for the possible factors $\mu$.\\
\\
{\bf Condition:} The factor $\mu$ divides (the squarefree part of) 
$D=c_0\psi_{01}^{(0)}+c_1\psi_{01}^{(1)}+c_2\psi_{01}^{(2)}$. \\
\\
A second (trivial) necessary condition for $\mu$ is given by the
property that the two quadratic equation\\  
\parbox{.05\textwidth}{(24)} 
\parbox{.95\textwidth}{
  \begin{alignat*}{3} 
  & \mu\sigma_0^2\ &&=\ \Psi_0(\eta_0,\eta_1)\ &&=\ \psi_{00}^{(0)}\eta_0^2 + \psi_{01}^{(0)}\eta_0\eta_1 
          + \psi_{11}^{(0)}\eta_1^2, \\ 
  & \mu\sigma_1^2\ &&=\ \Psi_1(\eta_0,\eta_1)\ &&=\ \psi_{00}^{(1)}\eta_0^2 + \psi_{01}^{(1)}\eta_0\eta_1 
          + \psi_{11}^{(1)}\eta_1^2
  \end{alignat*} } \\  
are both (individually) solvable (in the variables $\eta_0,\eta_1$ and $\sigma_0$ or $\sigma_1$, 
respectively. Apart from the above necessary conditions for the factor $\mu$, there seems to be 
no {\it a priori} criterion for selecting an appropriate factor $\mu$. We have to resort to a 
trial-and-error approach going through the possible candidates. We note in passing that the 
logarithmic height of $(Y_0,Y_1,Y_2)$ is about twice the logarithmic height of $(\eta_0,\eta_1)$.

\subsubsection{New quadratic form}

Having chosen a suitable factor $\mu$ according to 4.2.2, the original
problem is transformed to solving the system $(24)$. (Note that the orders of magnitude 
of the parameters $\eta_0$, $\eta_1$ and of $\sigma_0$, $\sigma_1$ are approximately equal.) 
\smallskip 

Let us rename the variables $(\eta_0,\eta_1,\sigma_0)$ to $(Z_0,Z_1,Z_2)$ and let us denote 
the first of the quadratic equations in $(24)$ by $Q_4(Z_0,Z_1,Z_2)=0$. Again, let us determine
a rational point $(z_0,z_1,z_2)$ on this quadric and let us parametrize $Q_4$ by using this point 
as a projection point. Then we obtain a representation 
$$Z_0\ =\ \Gamma_0(\rho_0,\rho_1), \qquad Z_1\ =\ \Gamma_1(\rho_0,\rho_1), \qquad 
  Z_2\ =\ \Gamma_2(\rho_0,\rho_1) \leqno{(25)}$$  
with quadratic forms $\Gamma_i$. Again, for any pair of parameters $(\rho_0,\rho_1)$ we get a 
solution for the first equation $Q_4(Z_0,Z_1,Z_2)=0$. Substituting the expressions for $Z_0$ 
and $Z_1$ into the right-hand side of the second equation yields 
$$\Psi_1\bigl(\Gamma_0(\rho_0,\rho_1),\Gamma_1(\rho_0,\rho_1)\bigr)\ =\ \mu\sigma_1^2. 
  \leqno{(26)}$$
This equation (which is of degree $4$ in the variables $\rho_0,\rho_1$) is the base for a 
final search loop. We note that the logarithmic height of $(Z_0,Z_1,Z_2)$ is about twice 
the logarithmic height of $(\rho_0,\rho_1)$.

\subsubsection{Algorithm}

The considerations above may be combined to the following strong form
of the algorithm.\\ 
\\
\eins\texttt{/{*} strong algorithm {*}/} \hfill\break 
\eins\texttt{INPUT: quadrics $Q_{1}$ and $Q_{2}$} \hfill\break 
\eins\texttt{BEGIN} \hfill\break
\eins\texttt{calculate a representation of the quadric intersection as in 4.2.1}\hfill\break
\eins\texttt{calculate a special point on $Q_1$ of the form $x=(0,x_1,x_2)$} \hfill\break 
\eins\texttt{calculate the parametrization of $Q_1$ projecting from $x$:} \hfill\break 
\zwei\texttt{$X_{i}=\Phi_{i}(\xi_{0},\xi_{1})$, $i=0,1,2$;}\hfill\break
\eins\texttt{substitute $X_i=\Phi_i(\xi_0,\xi_1)$, $i=0,1$, and plug into $Q_2$} \hfill\break 
\eins\texttt{calculate quadratic form $Q_3(Y_0,Y_1,Y_2)$ from this substitution} \hfill\break 
\eins\texttt{determine point $y=(y_0,y_1,y_2)$ on $Q_3$} \hfill\break 
\eins\texttt{calculate the parametrization of $Q_3$ projecting from $y$:} \hfill\break 
\zwei\texttt{$Y_i=\Psi_i(\eta_0,\eta_1)$, $i=0,1,2$} \hfill\break
\eins\texttt{determine the system of linear equations (22)} \hfill\break
\eins\texttt{determine a solution of this system} \hfill\break 
\eins\texttt{calculate the value $D$ from 4.2.2} \hfill\break 
\eins\texttt{determine the possible values for the coefficient $\mu$} \hfill\break 
\eins\texttt{choose an appropriate value for $\mu$} \hfill\break 
\eins\texttt{calculate the quadratic form $Q_4(Z_0,Z_1,Z_2)$ from 4.2.3} \hfill\break 
\eins\texttt{determine a point $z=(z_0,z_1,z_2)$ on $Q_4$} \hfill\break 
\eins\texttt{calculate the parametrization of $Q_4$ projecting from $z$:} \hfill\break 
\zwei\texttt{$Z_i=\Gamma_i(\rho_0,\rho_1)$, $i=0,1,2$} \hfill\break 
\eins\texttt{calculate the quartic form $\Psi_1(\Gamma_0(\rho_0,\rho_1),
    \Gamma_1(\rho_0,\rho_1))$} \hfill\break
\eins\texttt{loop over $(\rho_0,\rho_1)$} \hfill\break
\zwei\texttt{BEGIN} \hfill\break 
\drei\texttt{calculate the value val$=\Psi_1\bigl(\Gamma_0(\rho_0,\rho_1),\Gamma_1(\rho_0,\rho_1)\bigr)$}
  \hfill\break
\drei\texttt{check whether $\mu\cdot$val is a square} \hfill\break 
\drei\texttt{if yes then} \hfill\break 
\drei\texttt{BEGIN} \hfill\break 
\vier\texttt{/{*} Solution found {*}/} \hfill\break 
\vier\texttt{BREAK} \hfill\break 
\drei\texttt{END} \hfill\break 
\zwei\texttt{END} \hfill\break 
\eins\texttt{calculate values $Z_i=\Gamma_i(\rho_0,\rho_1)$, $i=0,1,2$} \hfill\break 
\eins\texttt{SET $\eta_i=Z_i$ , $i=0,1$} \hfill\break 
\eins\texttt{calculate values $Y_i=\Psi_i(\eta_0,\eta_1)$} \hfill\break 
\eins\texttt{calculate $\xi_i=\sqrt{\mu Y_i}$, $i=0,1$} \hfill\break 
\eins\texttt{calculate values $X_i=\Phi_i(\xi_0,\xi_1)$, $i=0,1,2$} \hfill\break 
\eins\texttt{if necessary, rename the variables as explained in 4.2.1} \hfill\break 
\eins\texttt{RETURN $(X_0,X_1,X_2,X_3)$} \hfill\break 
\eins\texttt{END}\\
\\
{\bf Remarks:} 
\begin{itemize} 
\item  As was the case with its weak form, the algorithm is not guaranteed to terminate. 
\item The logarithmic height of a solution $(x_0,x_1,x_2,x_3)$ found by the algorithm is about 
  four times the logarithmic height of the parameters $(\rho_0,\rho_1)$ of the final loop. Hence 
  the number of digits in the solutions which can be found by a search up to a certain height for 
  the parameters in the final search loop is about twice as large as compared to the weak form 
  of the algorithm. This improvement is counteracted by the fact that the additional condition 
  of 4.2.1 needs to be imposed on the given quadratic forms. Another disadvantage is the task of 
  choosing a factor $\mu$ according to 4.2.3, for which we are not aware of a guiding principle. 
\item As is the case with its weak form, the complexity of the algorithm depends quadratically 
  on the height of the parameters $(\rho_0,\rho_1)$. 
\end{itemize} 
\section{Examples}

This section is devoted to the discussion and presentation of explicit
examples to which the above algorithm is applied. In the first example 
we demonstrate the execution of the algorithm, exploiting the choices to 
be made for making the algorithm work and presenting the intermediate
data produced within the algorithm. The second subsection collects some 
series of examples which are very similar in their behaviour, while the 
third part treats examples with rank larger than one.

\subsection{An explicit example}

Let us consider the case of the $2\pi/3$-congruent number $n=142=2\cdot 71$. 
(Note that $71$ is a prime number congruent to $-1$ mod $8$ and hence is 
$2\pi/3$-congruent; cf. [24].) 
The elliptic curve associated with this problem is given by $y^2=x(x+142)
(x-3\cdot 142)$, and the equations for the corresponding concordant form problem 
are $W_0^2-3\cdot 142\,W_1^2=W_2^2$ and $W_0^2+142\,W_1^2=W_3^2$. From the 
equivalence classes determined by the 2-descent we get the class belonging to 
the triplet $(A,B,C)=(1,2,2)$, which is a good candidate for yielding a 
homogeneous space having a rational point. The corresponding equations are
$$x+142=\alpha^2,\qquad x=2\beta^2, \qquad x-3\cdot 142=2\gamma^2.
  \leqno{(27)}$$ 
These equations are equivalent to the following system of quadratic
equations
$$Q_1:\ 3X_0^2-8X_1^2+2X_2^2=0, \qquad Q_2:\ X_0^2-2X_1^2-142X_3^2=0. 
  \leqno{(28)}$$ 
Each of these equations individually possesses rational solutions; so the
corresponding homogeneous space has a chance to have a rational point. We see 
that $Q_1$ has the rational point $(0,1,2)$, so the critical condition of 
section 4.2 is fulfilled. Parametrizing $Q_1$ with this point gives 
$$X_0\ =\ 16\xi_0\xi_1, \qquad X_1\ =\ 8\xi_0^2+3\xi_1^2, \qquad 
  X_2\ =\ -16\xi_0^2 + 6\xi_1^2. \leqno{(29)}$$ 
Substituting the expressions for $X_0$ and $X_1$ in $Q_2$ and setting 
$Y_0=\xi_0^2$, $Y_1=\xi_1^2$, $Y_2=X_3$ gives the quadric 
$$Q_3:\ -64Y_0^{2}+80Y_0Y_1-9Y_1^2-71Y_2^2\ =\ 0.\leqno{(30)}$$ 
The next step is to find a point on $Q_3$. We find the point $(y_0,y_1,y_2) 
= (10,9,1)$. Parametrizing $Q_3$ using this point gives the following 
representations:\\
\parbox{.05\textwidth}{(31)} 
\parbox{.95\textwidth}{
  \begin{alignat*}{3} 
  & Y_0\ &&=\ \Psi_0(\eta_0,\eta_1)\ &&=\ -90\eta_0^2 + 81\eta_0\eta_1 - 20\eta_1^2, \\
  & Y_1\ &&=\ \Psi_1(\eta_0,\eta_1)\ &&=\ -719\eta_0^2 + 640\eta_0\eta_1 - 144\eta_1^2,  \\ 
  & Y_2\ &&=\ \Psi_2(\eta_0,\eta_1)\ &&=\ -9\eta_0^2 + 40\eta_0\eta_1 - 16\eta_1^2. 
  \end{alignat*} } \\
Now we have to find an appropriate coefficient $\mu$ which, according to the 
necessary condition from the linear system in 4.2.2, must divide $D=142$. The 
only candidates for $\mu$ which result in (individually) solvable quadratic 
equations $\Psi_0(\eta_0,\eta_1) = \mu\sigma_0^2$ and $\Psi_1(\eta_0,\eta_1) = 
\mu\sigma_1^2$ are $\mu=-1$ and $\mu=-71$. We do not have an {\it a priori} 
criterion which of these values to use for the following search, but must use 
trial and error. It turns out that $\mu = -1$ does not yield a solution, but 
$\mu=-71$ does. So let us consider this choice. According to the notations of 
section 4, we denote by $Q_4$ and $Q_5$ the quadratic equations $\Psi_0
(\eta_0,\eta_1)=\mu\sigma_0^2$ and $\Psi_1(\eta_0,\eta_1)=\mu\sigma_1^2$, 
respectively, where we rename the variables as follows: $Z_0=\eta_0$, $Z_1
=\eta_1$, $Z_2=\sigma_0$, $Z_3=\sigma_1$. Thus we arrive at the following equations:\\ 
\parbox{.05\textwidth}{(32)} 
\parbox{.95\textwidth}{
  \begin{alignat*}{3} 
  &Q_4:&&\quad -90Z_0^2 + 81Z_0Z_1 - 20Z_1^2+71Z_2^2\ &&=\ 0, \\
  &Q_5:&&\quad -719Z_0^2 + 640Z_0Z_1 - 144Z_1^2 + 71Z_3^2\ &&=\  0.
  \end{alignat*}} \\  
Again, we have to look for a point on $Q_4$. We find the point $(z_0,z_1,z_2)=(4,1,4)$ 
and use this point to parametrize $Q_4$. This yields \\
\parbox{.05\textwidth}{(33)} 
\parbox{.95\textwidth}{
  \begin{alignat*}{3} 
  & Z_0\ &&=\ \Gamma_0(\rho_0,\rho_1)\ &&=\ -5\rho_0^2 + 10\rho_0\rho_1 + 279\rho_1^2, \\ 
  & Z_1\ &&=\ \Gamma_1(\rho_0,\rho_1)\ &&=\ -19\rho_0^2 + 180\rho_0\rho_1 - 90\rho_1^2, \\ 
  & Z_2\ &&=\ \Gamma_2(\rho_0,\rho_1)\ &&=\ -5\rho_0^2 + 81\rho_0\rho_1 - 360\rho_1^2.  
  \end{alignat*}} \\  
Finally, we substitute the expressions for $Z_0$ and $Z_1$ into $Q_5$; 
this yields: 
$$Q_5:\ \begin{matrix} 
  \Psi_1(Z_0,Z_1) + 71Z_3^2\ =\ &-9\,159\,\rho_0^4 + 359\,260\,\rho_0^3\rho_1 
                - 5\,176\,610\rho_0^2\rho_1^2\phantom{\ =\ 0} \\ 
  &+ 32\,218\,380\,\rho_0\rho_1^3 - 73\,204\,479\,\rho_1^4 + 71\,Z_3^2\ =\ 0 
    \end{matrix} \leqno{(34)}$$ 
This last expression is the starting point for the final search loop. We look
for integers $(\rho_0,\rho_1)$ such that 
$$\begin{matrix}
  -71\,\Psi_1(\Gamma_0(\rho_0,\rho_1),\Gamma_1(\rho_0,\rho_1))\ =\ &-9\, 159\,\rho_0^4 
      + 359\,260\,\rho_0^3\rho_1 - 5\,176\,610\,\rho_0^2\rho_1^2 \\ 
  &+ 32\,218\,380\,\rho_0\rho_1^3 - 73\,204\,479\,\rho_1^4 
   \end{matrix} \leqno{(35)}$$ 
is a perfect square. Having found such a pair of integers we can recover
a solution of the original equations. Using $(\rho_0,\rho_1)=(20,3)$, the 
right-hand side of (35) takes the value $-10633599 = -71\cdot 387^2$, so 
we have found a solution. With these values we find that \\
\parbox{.05\textwidth}{(36)}
\parbox{.95\textwidth}{
  \begin{alignat*}{2} 
  & (Z_0,Z_1,Z_2,Z_3)\ &&=\ (1111, 2390, -380, 387) \\ 
  & (Y_0,Y_1,Y_2)    \ &&=\ (-10252400, -10633599, 3709111) \\ 
  & (X_0,X_1,X_2,X_3)\ &&=\ (2352960, 1604507, -1411786, -52241)
  \end{alignat*}} \\
which results in the point 
$$(x,y)\ =\ \left( \frac{5148885426098}{2729122081}, 
   \frac{10659946547134851840}{142572066633521}\right) \leqno{(37)}$$
on the curve $y^2=x(x+142)(x-426)$. The corresponding solution 
$(W_0,W_1,W_2,W_3)$ of the concordant form problem is given by  
$$\left[\begin{matrix} W_0\\ W_1\\ W_2\\ W_3\end{matrix}\right]\ =
  \ \left[\begin{matrix} -1685098252492020382767601 \\
                         \phantom{-}69610783446108974371680 \\ 
                         -880513748494434998396401 \\ 
                         -1878201269026558326761999 \end{matrix}\right]  
  \leqno{(38)}$$ 
We note that the signs in this solution are not important; they arise from the 
isomorphism used between the elliptic curve and the quadric intersection defining 
the concordant form problem. Other distributions of signs correspond to other 
points on the elliptic curve; these other points are determined by adding any of 
the 2-torsion points of $E_{1,3,142}(\mathbb{Q})$ to the above solution. They are 
given by the following three points, together with their negatives (in the sense 
of the elliptic curve addition), which are obtained by reversing the sign of the 
$y$-component: \\
\parbox{.05\textwidth}{(39)} 
\parbox{.95\textwidth}{
 \begin{eqnarray*} 
 (x_1,y_1) &=& \left( \frac{-82545026461926}{2574442713049},
               \frac{5248834080776243516160}{4130711354186111843}\right), \\
 (x_2,y_2) &=& \left( \frac{294814405555200}{498284927449},
               \frac{2982672665844557232960}{351735842291756957}\right), \\
 (x_3,y_3) &=& \left( \frac{-35378229848879}{346026297600},
               \frac{298269379294025686631}{203546509300224000}\right).
 \end{eqnarray*}} \\  
{\bf Remarks.} During the execution of the algorithm one has to make some choices, 
most of which are not very critical and do not affect the working of the algorithm. 
Examples of such choices are those of the points on the quadrics $Q_1$, $Q_3$ and 
$Q_4$, which are used for the parametrizations of these quadrics. Also, it is 
irrelevant whether one uses the quadric $Q_4$ for parametrization and substitution 
into $Q_5$ or whether one exchanges the roles of $Q_4$ and $Q_5$. A more delicate 
task is the choice of the parameter $\mu$, for which we have no good suggestions. 
Note that using different parametrizations of the intermediate quadrics may affect 
the possible choices for this parameter $\mu$.

\subsection{Some series of examples with similar behaviour}

\subsubsection{Congruent prime numbers $k\equiv5$ mod $8$}

Let $k\in\mathbb{N}$ be a prime number satisfying $k\equiv5$ mod $8$. Then
$k$ is a congruent number, i.e., it is the area of a right triangle with 
rational sides. This is tantamount to saying that there are nontrivial rational 
points on the elliptic curve $y^2=x(x+k)(x-k)$, and the corresponding concordant 
form problem is given by the equations $W_0^2-kW_1^2=W_2^2$ and $W_0^2+kW_1^2=W_3^2$. 
The Mordell-Weil-rank of these curves is one, and rational points on these curves 
can be found by examining the 2-descent with parameters $(A,B,C)=(1,-1,-1)$. The 
initial quadrics for the algorithm are given by
$$Q_1:\ 2X_0^2+X_1^2-X_2^2=0, \qquad Q_2:\ X_0^2+X_1^2-kX_3^{2}=0.\leqno{(40)}$$ 
The point $(0,1,1)$ is on $Q_1$ and can be used for the first parametrization.
Appendix 1 provides a table which contains the results found for all prime 
numbers $k\equiv 5$ mod $8$ up to $613$ in terms of the solutions $(W_0,W_1,W_2,W_3)$ 
to the concordant form problem (with all the components being positive). The table 
shows (in the last column) also the logarithmic heights of the solutions found 
(with respect to the logarithm with base $10$, i.e., the maximum number of
decimal places of the solution). The following diagram shows the logarithmic 
heights depending on the prime numbers $k$.

\includegraphics[width=12cm]{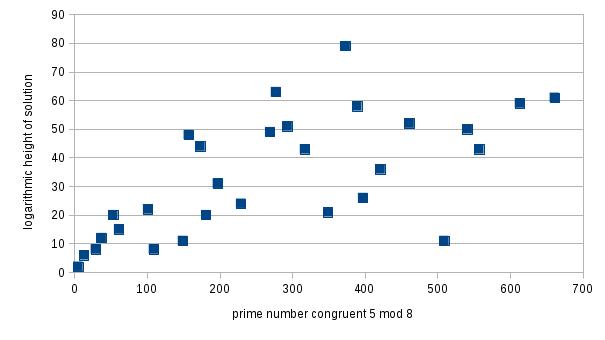}

\subsubsection{Diagrams of some other series}

$\phantom{con}$(a) The following diagram shows the heights of the solutions (in terms  
of the solutions $(W_0,W_1,W_2,W_3)$ of the concordant form
problem) depending on the coefficients $k$, where $k$ is a prime
number satisfying $k\equiv 7$ mod $8$.

\includegraphics[width=12cm]{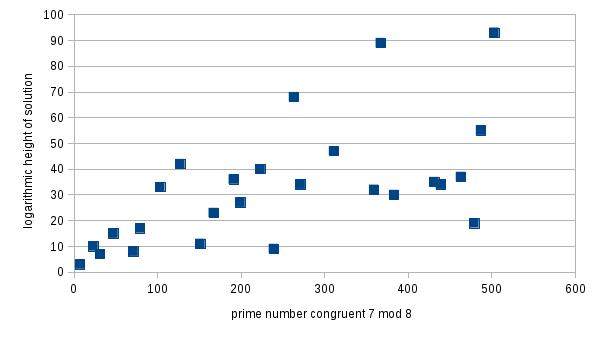}\\

(b) The following diagram shows the logarithmic heights of the solutions
to the congruent number problem for the numbers $k=2\ell$ where $\ell$ is a  
prime number satisfying $\ell\equiv 7$ mod $8$.

\includegraphics[width=12cm]{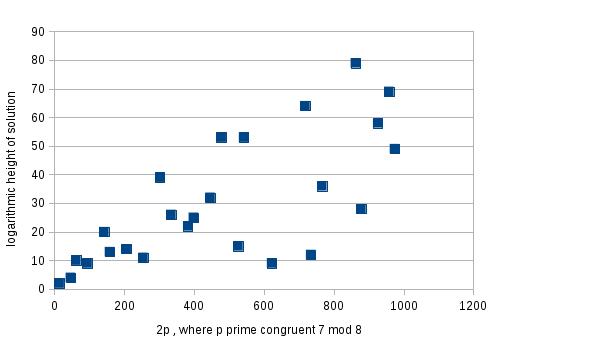}\\

(c) Finally, let us consider the $2\pi/3$-congruent number problem. This 
corresponds to the elliptic curves $E_{1,3,k}$ given by $y^2=x(x+k)(x-3k)$.
From [24] we know that prime numbers $k\equiv 5$ mod $24$ are $2\pi/3$-congruent. 
The following diagram shows again the logarithmic heights of solutions in terms 
of the solutions to the corresponding concordant form problem.

\includegraphics[width=12cm]{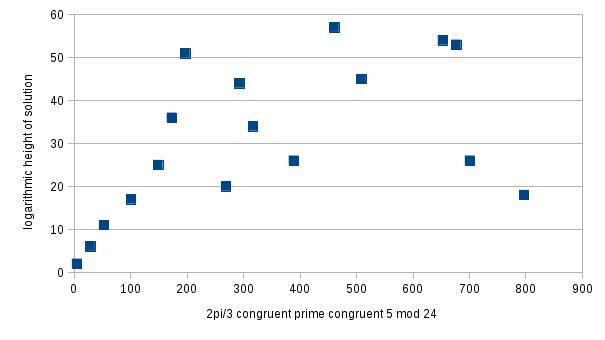}

\subsection{Special examples}

Let $r=2\ell$ where $\ell\in\mathbb{N}$ is a prime number satisfying $\ell\equiv 7$
mod $96$. As in section 3, we consider the $2\pi/3$-congruent number problem for 
these numbers. The solution to this problem corresponds to determining the rational 
points on the elliptic curve $y^2=x(x+2\ell)(x-6\ell)$, and the corresponding 
concordant form problem is given by the equations $W_0^2+2\ell W_1^2 = W_2^2$ and 
$W_0^2-6\ell W_1^2 = W_3^2$. According to the considerations in section 3, we have 
a good chance to find rational solutions to this problem by examining the 
homogeneous spaces belonging to the 2-descent parameter sets $(A,B,C)=(1,2,2)$,
$(2,-3,-6)$ and $(2,-6,-3)$, where the solutions to the third parameter set are 
obtained as the sums of the solutions of the first two parameter sets in terms 
of the group structure of the elliptic curve. Appendix 2 shows the results obtained 
for the three examples $\ell=7$, $\ell=103$ and $\ell=199$. For the concordant form 
we only exhibit the solution with positive coefficients; in addition, we list 
the four points on the elliptic curve with positive $y$-coordinates which 
correspond to the concordant form solutions obtained by changing the signs of 
the coefficients (which means adding 2-torsion points on the elliptic curve).  
\section{Concluding remarks and open questions}

The algorithm described in section 4 provides us with a tool to find 
explicit solutions to elliptic curve equations of a special form.
This algorithm improves a simpler strategy which is probably well
known to the experts and which is based on a classical 2-descent procedure
together with a parametrization scheme for quadratic forms. As was shown 
by way of various examples, the algorithm works quite well in many situations. 
However, there are certain points which are not fully understood, and these 
open questions are addressed in this final section. 

\subsection{Choices within the algorithm}

The algorithm applies to elliptic curves with full 2-torsion $\mathbb{Z}/2\mathbb{Z}
\times\mathbb{Z}/2\mathbb{Z}$. It requires the validity of a certain condition which 
seems to be of a purely technical nature; namely, we assumed that the homogeneous 
space leading to a solution of the elliptic curve equation can be given by two 
diagonal quadrics with separated variables in the form \\
\parbox{.05\textwidth}{(41)} 
\parbox{.95\textwidth}{
  \begin{alignat*}{3} 
  & Q_1: && \quad a_{00}X_0^2 + a_{11}X_1^2 + a_{22}X_2^2\ &&=\ 0 \\
  & Q_2: && \quad b_{00}X_0^2 + b_{11}X_1^2 + b_{33}X_3^2\ &&=\ 0 
  \end{alignat*}} \\ 
such that at least one of the equations (without loss of generality
the first one) has a solution of the form $(0,x_1,x_2)$ or $(x_0,0,x_2)$. 
However, this condition is not always fulfilled, as is shown by the following 
example. \smallskip 

Let $k$ be a prime number congruent to $23$ modulo $24$. Then, assuming 
the validity of the Birch-Swinnerton-Dyer Conjecture, the elliptic curve 
$E_{1,3,k}$ has positive rank $\bigl($cf. [24], section 3, expectation 
(e1)$\bigr)$; in fact, we expect the rank of this curve to be one. Let us consider  
the special case $k=23$. Then by the 2-descent procedure we find the triplet 
$(A,B,C)=(2,3,6)$ to determine a homogeneous space having a rational solution. 
The four quadrics which are determined by this triplet according to 4.2.1 
are given by\\
\parbox{.05\textwidth}{(42)} 
\parbox{.95\textwidth}{
  \begin{alignat*}{3} 
  & 0\ &&=\ X_0^2 - 2X_1^2 + X_2^2 &&\\
  & 0\ &&=\ 2X_0^2 - 3X_1^2 - 23 X_3^2\ &&=\ 2X_0^2 - 3X_1^2 - kX_3^2 \\ 
  & 0\ &&=\ X_0^2 - 3X_2^2 - 46 X_3^2\ &&=\ X_0^2 - 3 X_2^2 - 2kX_3^2 \\ 
  & 0\ &&=\ X_1^2 - 2X_2^2 - 23 X_3^2\ &&=\ X_1^2 - 2X_2^2 - kX_3^2
  \end{alignat*}} \\ 
and it is easy to see that none of these equations has an integer 
solution with one component being $0$. However, $(7,5,1,1)$ is a
solution for this system of quadratic equations, and this solution
determines for example the point $(75,210)$ on the elliptic curve 
$E_{1,3,23}$. The same equations yield solutions for other values of $k$, 
for example $k = 47, 71, 167, 191, 239, 263, 311, 359, 383, 431, 479, 503, 599$. 
\smallskip 

Note that for any of the triplets $(A',B',C')$ which are equivalent to $(2,3,6)$ 
with respect to the 2-torsion points, the critical condition is not satisfied 
either; hence we cannot apply the strong form of the algorithm to find rational 
points on the curves of the form $E_{1,3,k}$ with prime numbers $k\equiv 23$ 
modulo $24$. \\ 
\\ 
{\bf Question:} Is there a geometric or arithmetic interpretation of 
the above condition? \\

During the execution of the algorithm we have to make some choices, some of which 
cause no problems at all whereas others need to be treated carefully.\\ 

1. We have to choose an appropriate homogeneous space leading to a
solution. There are at least the homogeneous spaces being $P$-equivalent
in the sense of section 2, where $P$ is the set of 2-torsion points
on the elliptic curve. The choice of different spaces should yield
different solutions which are in the same residue class with respect
to the 2-torsion subgroup.\\ 
 
2. In several situations we have to choose a point on a quadric for
defining a parametrization of that quadric. Since the parametrizations
always determine all rational points on that quadric, the choice of
the point is irrelevant. However, the solutions found by choosing
different points may arise in different succession.\\

3. At a sensitive point of the algorithm, we have to choose a factor
$\mu$ which determines the ``right'' representative of a point in
projective space to have square coefficients. We have determined a
finite set of possibilities for this choice, and we can exclude many of
the potential candidates by simple arguments. However, for the remaining
possibilities we have no guidelines which of the candidates may give rise 
to a solution, and hence resort to a trial-and-error approach.\\
\\
{\bf Question:} Is there a geometric or arithmetic interpretation of this
factor which may lead to a better way of choosing this factor? 

\subsection{Qualities of the algorithm}

Summarizing the observations made in section 4, we see that a search within 
the algorithm in which the parameters $\rho_0$ and $\rho_1$ of the  central loop 
are considered up to logarithmic height $n$, the solutions will have logarithmic 
height of about $12n$. Hence a search with parameters up to around $10\,000$ is 
expected to find solutions with about 50 decimal places. Such a search can be 
performed with a usual personal computer in a few minutes. However, as mentioned 
in section 4, the complexity of the algorithm grows quadratically with the parameters 
$(\rho_0,\rho_1)$. Hence a search in a parameter range of up to $100\, 000$, which 
could provide solution with about $70$ decimal places, would already take several 
hours.\smallskip 

In spite of this complexity defect, the examples of section 5 show that we can 
generate enough data to obtain interesting information on the behaviour of the 
smallest solutions found by the algorithm, depending on the coefficients of the 
original elliptic curves. Of course the data are too sparse to lead to serious 
conjectures. Nevertheless, the diagrams in section 5 suggest that for the families
considered there may be a linear correlation between the parameters
defining the curves and the logarithmic height of the smallest solutions.
We are far from being able to formulate such a suspected correlation
in more precise terms.

\subsection{Structural questions}

The central point in both the weak as well as the strong version of the 
algorithm is the consideration of a suitable homogeneous space $Q$ over 
the given elliptic curve $E$. Such a homogeneous space is geometrically 
(i.e. over an algebraic closure of $\mathbb{Q}$) isomorphic to $E$. It is 
trivial (in the sense of the Weil-Chatelet group) if and only if it contains 
a rational point, and in this situation it is isomorphic over $\mathbb{Q}$ to 
some rationally defined elliptic curve $E_{Q}$. Via the well known construction 
of Nagell (cf. [15]) we can define a biregular mapping from $Q$ to $E_Q$ given 
by a Weierstra\ss{} equation. \\
\\
{\bf Question:} What can be said about this elliptic curve?\\
\\
We observe that if rank$(E_{Q}(\mathbb{Q}))>0$ we obtain an infinite
series of rational points on $E$ by considering the composition of
an isomorphism of $Q$ to $E_{Q}$ and the mapping from $Q$ to $E$ given in 
section 2. Obviously, independent rational points on $E_{Q}(\mathbb{Q})$
will determine independent points on $E(\mathbb{Q})$. So we certainly have 
$\hbox{rk}(E_{Q}(\mathbb{Q}))\leq\hbox{rk}(E(\mathbb{Q}))$. In the 
examples of section 5 two cases occurred. 
\begin{itemize} 
\item There were examples with $\hbox{rk}(E(\mathbb{Q}))=1$, and hence for the 
  choosen homogeneous space to yield a solution of infinite order we also had 
  necessarily $\hbox{rk}(E_{Q}(\mathbb{Q}))=1$.
\item There were examples with $\hbox{rk}(E(\mathbb{Q}))=2$, and in these examples 
  invariably all the homogeneous spaces $Q$ considered for finding solutions on 
  $E(\mathbb{Q})$ were such that $\hbox{rk}(E_{Q}(\mathbb{Q}))=1$. Hence to find 
  independent solutions we had to look for two different (and non-equivalent with 
  respect to 2-torsion points) homogeneous spaces.
\end{itemize} 
{\bf Question:} Are there examples in which $\hbox{rk}(E_{Q}(\mathbb{Q}))>1$?\\
\\
Let us again consider the case in which $\hbox{rk}(E(\mathbb{Q}))\geq 2$, and
let $Q$ and $Q'$ be two independent homogeneous spaces yielding independent 
rational solutions on $E(\mathbb{Q})$. Then the associated rationally
defined elliptic curves $E_{Q}$ and $E_{Q'}$ are twists of one another (and 
twists of the original curve $E$ as well). \\
\\
{\bf Question:} What can be said about the connection between these two elliptic
curves? 

\subsection{Further outlook}

The algorithm above was developed for the class of elliptic curves
$E$ with full 2-torsion, which means that the torsion subgroup of $E$ 
contains $\mathbb{Z}/2\mathbb{Z}\times\mathbb{Z}/2\mathbb{Z}$. Equivalently, 
the curves considered can be given in affine form by a Weierstra\ss{}
equation with split polynomial in the form $y^2=(x-e_1)(x-e_2)(x-e_3)$ 
where $e_1,e_2,e_3\in\mathbb{Q}$ are pairwise different. It would be of 
interest to develop algorithms for other classes of elliptic curves with 
positive rank. Note that the 2-descent procedure may always be formulated 
in the same way as in section 2, but over some number field instead of 
$\mathbb{Q}$. One could try to develop an algorithm over such a number field 
and then, {\it a fortiori}, try to extract those solutions which are actually 
rational. Alternatively, one could try to reformulate the 2-descent procedure 
(at least in special cases) such that a rationally defined algorithm can be 
developed. Both strategies seem to be not entirely trivial. \\
\vskip 2 true cm 
\noindent 
Hagen Knaf, Karlheinz Spindler \hfill\break
Hochschule RheinMain, Germany \hfill\break
Applied Mathematics \hfill\break
{\tt hagen.knaf@hs-rm.de}, {\tt karlheinz.spindler@hs-rm.de} 
\par\bigskip\noindent 
Erich Selder \hfill\break
Frankfurt University of Applied Sciences, Germany \hfill\break
Computer Science and Engineering \hfill\break 
{\tt e\_selder@fb2.fra-uas.de} 

\vfill\eject 

\section{References} 

\begin{enumerate} 

\item Andrew Bremner, John William Scott Cassels, 
\emph{On the equation $Y^{2}=X(X^{2}+p)$},
Mathematics of Computation, Vol. 42, No. 165, 1984, pp. 257--264.

\item John William Scott Cassels, 
\emph{Lectures on Elliptic Curves}, 
London Mathematical Society Students Texts vol. 24, London 1991. 

\item Leonhard Euler, 
\emph{De binis formulis speciei xx+myy et xx+nyy inter se concordibus et discordibus}, 
Mem. Acad. Sci. St.-Petersbourg 1780 (Opera Omnia: Ser. 1, Vol. 5, pp. 406--413).

\item Masahiko Fujiwara,
\emph{$\theta$-congruent numbers},
in: K. Gy\"o{}ry et al. (eds.), \emph{Number theory},   
de Gruyter, Berlin 1998, pp. 235--241.

\item Masahiko Fujiwara,
\emph{Some properties of $\theta$-congruent numbers},
Natural Science Report, Ochanomizu University,   
vol. 52, no. 2, 2001.

\item Ludwig Holzer, 
\emph{Minimal solutions of diophantine equations}, 
Canad. J. Math. 2, 1950, pp. 238--244.

\item Makiko Kan,
\emph{$\theta$-congruent numbers and elliptic curves},  
Acta Arithmetica 94 (2), 2000, pp. 153--160.

\item Anthony W. Knapp, 
\emph{Elliptic Curves}, Mathematical Notes 40, 
Princeton University Press 1992.

\item Neal Koblitz,
\emph{Introduction to Elliptic Curves and Modular Forms}, 
Springer, New York/Berlin/Heidelberg 1993.

\item Adrien-Marie Le Gendre, 
\emph{Recherches d'Analyse ind\'etermin\'ee}, 
Histoire de l'Aca\-d\'emie Royale des Sciences 1785, 
pp. 465--559. 

\item Elisabeth Lutz, 
\emph{Sur l'\'equation $y^2=x^3-Ax-B$ dans les corps $p$-adiques}, 
J. Reine Angew. Mathematik 177 (1937), pp. 237--247. 

\item Barry Charles Mazur, 
\emph{Modular curves and the Eisenstein ideal}, 
Publications math\'e{}matiques de l'I.H.E.S 47 (2), 1977, pp. 33--186.

\item Louis Joel Mordell, 
{\it On the Rational Solutions of the Indeterminate Equations of the 
Third and Fourth Degrees}, 
Proc. Cambridge Phil. Soc. XXI, 1922, pp. 179--192. 

\item Louis Joel Mordell, 
\emph{On the Magnitude of the Integer Solutions of the Equation $ax^2+by^2+cz^2=0$}, 
J. Number Theory 1, 1969, pp. 1--3. 

\item Trygve Nagell, 
\emph{Sur les propri\'et\'es arithm\'etiques des 
cubiques planes du premier genre}, Acta Mathematica 52 (1928), pp. 93--126.

\item Trygve Nagell, 
\emph{Solution de quelques probl\`emes dans la th\'eorie arithm\'etique des cubiques 
planes du premier genre}, Wid. Akad. Skrifter I, (1), 1935.  

\item Ken Ono, 
\emph{Euler's Concordant Forms}, 
Acta arithmetica LXXVIII (2), 1996, pp. 101--123.   

\item Takashi Ono, 
\emph{Variations on a Theme of Euler},  
Plenum Press, New York and London 1994.

\item Erich Selder, Karlheinz Spindler, 
\emph{On $\theta$-congruent numbers, rational squares in arithmetic progressions, 
concordant forms and elliptic curves}, 
Mathematics 3(1), 2015, pp. 2--15. 

\item Joseph H. Silverman, 
\emph{The Arithmetic of Elliptic Curves}, 
Springer 2009. 

\item Joseph H. Silverman, J. Tate, 
\emph{Rational Points on Elliptic Curves}, 
Springer, New York 1992. 

\item Andr\'e Weil: 
\emph{Sur un th\'eor\`eme de Mordell}, 
Bull. Sci. Math. 2 (54), 1930, pp. 182--191. 

\item Kenneth S. Williams, 
\emph{On the Size of a Solution of Legendre's Equation}, 
Utilitas Mathematica 34, 1988, pp. 65--72.

\item Shin-ichi Yoshida, 
\emph{Some Variants of the Congruent Number Problem I}, 
Kyushu J. Math. 55 (2001), pp. 387--404.  

\item Don Bernard Zagier, 
\emph{Elliptische Kurven: Fortschritte und Anwendungen}, 
JBer. DMV 92 (1990), pp. 58--76. 

\end{enumerate}
\vfill\eject

\def\hgt{\vphantom{\frac{\int^b}{\int_a}}}
\noindent 
Table 1: Solutions to the system $W_0^2-kW_1^2=W_2^2$, $W_0^2+kW_1^2=W_3^2$ \\ 
where $k$ is prime with $k\equiv 5$ mod $8$ \\
\begin{tabular}{|l|l|l|}
\hline 
$k$ & $W_{0},W_{1},W_{2},W_{3}$ & log\tabularnewline
 &  & hgt\tabularnewline
\hline 
\hline 
$5$ & $W_{0}=41$ & $2$\tabularnewline
 & $W_{1}=12$ & \tabularnewline
 & $W_{2}=31$ & \tabularnewline
 & $W_{3}=49$ & \tabularnewline
\hline 
$13$ & $W_{0}=106921$ & $6$\tabularnewline
 & $W_{1}=19380$ & \tabularnewline
 & $W_{2}=80929$ & \tabularnewline
 & $W_{3}=127729$ & \tabularnewline
\hline 
$29$ & $W_{0}=48029801$ & $8$\tabularnewline
 & $W_{1}=180180$ & \tabularnewline
 & $W_{2}=48019999$ & \tabularnewline
 & $W_{3}=48039601$ & \tabularnewline
\hline 
$37$ & $W_{0}=605170417321$ & $12$\tabularnewline
 & $W_{1}=9475102140$ & \tabularnewline
 & $W_{2}=602419674529$ & \tabularnewline
 & $W_{3}=607908713329$ & \tabularnewline
\hline 
$53$ & $W_{0}=4850493897329785961$ & $19$\tabularnewline
 & $W_{1}=595711308569957580$ & \tabularnewline
 & $W_{2}=2172343665411286111$ & \tabularnewline
 & $W_{3}=6506573990620136689$ & \tabularnewline
\hline 
$61$ & $W_{0}=250510625883241$ & $15$\tabularnewline
 & $W_{1}=18295510698660$ & \tabularnewline
 & $W_{2}=205760310228191$ & \tabularnewline
 & $W_{3}=288398755364209$ & \tabularnewline
\hline 
$101$ & $W_{0}=2015242462949760001961$ & $22$\tabularnewline
 & $W_{1}=118171431852779451900$ & \tabularnewline
 & $W_{2}=1628124370727269996961$ & \tabularnewline
 & $W_{3}=2339148435306225006961$ & \tabularnewline
\hline 
$109$ & $W_{0}=10537321$ & $8$\tabularnewline
 & $W_{1}=872340$ & \tabularnewline
 & $W_{2}=5299871$ & \tabularnewline
 & $W_{3}=13927729$ & \tabularnewline
\hline 
$149$ & $W_{0}=11880808361$ & $11$\tabularnewline
 & $W_{1}=879612300$ & \tabularnewline
 & $W_{2}=5086222111$ & \tabularnewline
 & $W_{3}=16013667889$ & \tabularnewline
\hline 
$157$ & $W_{0}=224403517704336969924557513090674863160948472041$ & $48$\tabularnewline
 & $W_{1}=17824664537857719176051070357934327140032961660$ & \tabularnewline
 & $W_{2}=21796977171070247104112455266586147721935979809$ & \tabularnewline
 & $W_{3}=316605068345983991287469841722668300352741098609$ & \tabularnewline
\hline 
\end{tabular}\\

\begin{tabular}{|l|l|l|}
\hline 
$k$ & $W_{0},W_{1},W_{2},W_{3}$ & log\tabularnewline
 &  & hgt\tabularnewline
\hline 
\hline 
$173$ & $W_{0}=11389552969201600543101928087171460571651881$ & $44$\tabularnewline
 & $W_{1}=151819892495256080406058239068697733204020$ & \tabularnewline
 & $W_{2}=11213134773123931932373766469330799882824031$ & \tabularnewline
 & $W_{3}=11563279908237839493160911313667068050342769$ & \tabularnewline
\hline 
$181$ & $W_{0}=10940671490772286441$ & $20$\tabularnewline
 & $W_{1}=812534430489915900$ & \tabularnewline
 & $W_{2}=447084261166681441$ & \tabularnewline
 & $W_{3}=15465985290352891441$ & \tabularnewline
\hline 
$197$ & $W_{0}=3976155246560604347409241506281$ & $31$\tabularnewline
 & $W_{1}=128879379273797845692300739620$ & \tabularnewline
 & $W_{2}=3540856019037985665622394486369$ & \tabularnewline
 & $W_{3}=4368290253857372604620867723569$ & \tabularnewline
\hline 
$229$ & $W_{0}=764646440211958998267241$ & $24$\tabularnewline
 & $W_{1}=9404506457489780613180$ & \tabularnewline
 & $W_{2}=751285786287393798649441$ & \tabularnewline
 & $W_{3}=777777618847556210645041$ & \tabularnewline
\hline 
$269$ & $W_{0}=3895373414239011964782976279255856376333539432681$ & $49$\tabularnewline
 & $W_{1}=27965347900755720997936131300398362642863770100$ & \tabularnewline
 & $W_{2}=3868276064680043715910003459047778706397631593569$ & \tabularnewline
 & $W_{3}=3922283564474102177246930010978679686556478603569$ & \tabularnewline
\hline 
$277$ & {\footnotesize{}$W_{0}=225651876701966818406248027783418906721100922839903398228891241$} & $63$\tabularnewline
 & {\footnotesize{}$W_{1}=7177227596170451913324498105378376615737197106454374393982740$} & \tabularnewline
 & {\footnotesize{}$W_{2}=191441323585574742208871474928771109013368200999037188254008609$} & \tabularnewline
 & {\footnotesize{}$W_{3}=255318934946162061510445203213880843051935526236604388007831409$} & \tabularnewline
\hline 
$293$ & $W_{0}=464650359520278159096671986562151812257229902698281$ & $51$\tabularnewline
 & $W_{1}=9525532939264666216445930388515870770466775878580$ & \tabularnewline
 & $W_{2}=435102716279337902002156546716850448714905366022431$ & \tabularnewline
 & $W_{3}=492428207448547244972204796522699373235322337048369$ & \tabularnewline
\hline 
$317$ & $W_{0}=7704952068030240987029060443439470576691561$ & $43$\tabularnewline
 & $W_{1}=273033470936425799912142469450375693280340$ & \tabularnewline
 & $W_{2}=5977859131736779268748519442377486506296289$ & \tabularnewline
 & $W_{3}=9110311352659588118420210652872068714343089$ & \tabularnewline
\hline 
$349$ & $W_{0}=543117687145297245481$ & $21$\tabularnewline
 & $W_{1}=24479594709742323420$ & \tabularnewline
 & $W_{2}=292981872551143852319$ & \tabularnewline
 & $W_{3}=710010751000672703281$ & \tabularnewline
\hline 
$373$ & {\scriptsize{}$W_{0}=6464736286838262275566375140640125524476830394378258160144359151221846588162921$} & $79$\tabularnewline
 & {\scriptsize{}$W_{1}=214402886988423616335778394508029972671920911384749815755228436417174376951980$} & \tabularnewline
 & {\scriptsize{}$W_{2}=4964526988887992094202607668810309975770378526931158358479760499172740751760929$} & \tabularnewline
 & {\scriptsize{}$W_{3}=7677180621382399924131415436519959747090354653821331133153517438341892919535729$} & \tabularnewline
\hline 
\end{tabular}\\

\begin{tabular}{|l|l|l|}
\hline 
$k$ & $W_{0},W_{1},W_{2},W_{3}$ & log\tabularnewline
 &  & hgt\tabularnewline
\hline 
\hline 
$389$ & {\footnotesize{}$W_{0}=7091795623967975164665712219283669343892955896357711001321$} & $58$\tabularnewline
 & {\footnotesize{}$W_{1}=178471843490509327250771615016308845953260187406885241100$} & \tabularnewline
 & {\footnotesize{}$W_{2}=6156546092792523632766658624365628021678212406860376224929$} & \tabularnewline
 & {\footnotesize{}$W_{3}=7917327235348034759843187663340785352885655138147390834929$} & \tabularnewline
\hline 
$397$ & $W_{0}=40610678141909645597145961$ & $26$\tabularnewline
 & $W_{1}=897770616925261772023980$ & \tabularnewline
 & $W_{2}=36458857951695016208049311$ & \tabularnewline
 & $W_{3}=44375737009650677093379889$ & \tabularnewline
\hline 
$421$ & $W_{0}=206116218357279640098356283784343401$ & $36$\tabularnewline
 & $W_{1}=5615337183197656507592648081062140$ & \tabularnewline
 & $W_{2}=170906168853566716230403514501943649$ & \tabularnewline
 & $W_{3}=236133166640367788597694235239072049$ & \tabularnewline
\hline 
$461$ & $W_{0}=3891001511194439641326936071293799433960980792636201$ & $52$\tabularnewline
 & $W_{1}=141603906393919705341026008387612184428936054638300$ & \tabularnewline
 & $W_{2}=2428183393618185294158471720281132279464765769008799$ & \tabularnewline
 & $W_{3}=4937986525618685709148575673193534971677351329281201$ & \tabularnewline
\hline 
$509$ & $W_{0}=8234822441$ & $11$\tabularnewline
 & $W_{1}=358112820$ & \tabularnewline
 & $W_{2}=1592388641$ & \tabularnewline
 & $W_{3}=11536416241$ & \tabularnewline
\hline 
$541$ & $W_{0}=20712649137553815516771958538277092457342029080681$ & $50$\tabularnewline
 & $W_{1}=692383596502537714323160801078335881010146641980$ & \tabularnewline
 & $W_{2}=13025402685121753757242473560733781819142094311519$ & \tabularnewline
 & $W_{3}=26236740527002218360424469656173951317763934832881$ & \tabularnewline
\hline 
$557$ & $W_{0}=5499709076648565793208509282424464890877481$ & $43$\tabularnewline
 & $W_{1}=73405451625480094969934137800528933306420$ & \tabularnewline
 & $W_{2}=5219720607933421868366276874910516409140831$ & \tabularnewline
 & $W_{3}=5766117986189355551406701674738483083808369$ & \tabularnewline
\hline 
$613$ & {\footnotesize{}$W_{0}=18030067140713632672003110416838548155838466663510382251241$} & $59$\tabularnewline
 & {\footnotesize{}$W_{1}=184995021722032435269683407294072813709395011761171749380$} & \tabularnewline
 & {\footnotesize{}$W_{2}=17438592982424790791434179315233932893903828538097188624609$} & \tabularnewline
 & {\footnotesize{}$W_{3}=18602744877856272407302202226070441326874205450445721695409$} & \tabularnewline
\hline 
\end{tabular}

\pagebreak{}
\noindent
Table 2: Solutions to the system $W_0^2+2\ell W_1^2=W_2^2$, $W_0^2-6\ell W_1^2=W_3^2$ \\ 
and to the equation $y^2=x(x+2\ell)(x-6\ell)$ where $\ell\in\{ 7,\, 103,\, 199\}$ \\
\begin{tabular}{|l|l|l|l|}
\hline 
$l$ & $(A,B,C)$ & $W_{0},W_{1},W_{2},W_{3}$ & $(x,y)$\tabularnewline
\hline 
\hline 
$7$ & $(1,2,2)$ & $W_{0}=193$ & $(50,160)$\tabularnewline
 &  & $W_{1}=20$ & $(-294/25,-4705/125)$\tabularnewline
 &  & $W_{2}=207$ & $(336,5880)$\tabularnewline
 &  & $W_{3}=143$ & $(-7/4,245/8)$\tabularnewline
\hline 
$7$ & $(2,-3,-6)$ & $W_{0}=61$ & $(-12,36)$\tabularnewline
 &  & $W_{1}=6$ & $(49,147)$\tabularnewline
 &  & $W_{2}=65$ & $(378,7056)$\tabularnewline
 &  & $W_{3}=47$ & $(-14/9,784/27)$\tabularnewline
\hline 
$7$ & $(2,-6,-3)$ & $W_{0}=13$ & $(-6,48)$\tabularnewline
 &  & $W_{1}=2$ & $(98,784)$\tabularnewline
 &  & $W_{2}=15$ & $(84,588)$\tabularnewline
 &  & $W_{3}=1$ & $(-7,49)$\tabularnewline
\hline 
$103$ & $(1,2,2)$ & $W_{0}=14497255873$ & $\hgt(\frac{-3757543}{33124},\frac{16707000155}{6028568})$\tabularnewline
 &  & $W_{1}=573225380$ & $\hgt(\frac{40941264}{36481},\frac{191037111720}{6967871})$\tabularnewline
 &  & $W_{2}=16669115727$ & $(\hgt\frac{470450}{289},\frac{269753120}{4913})$\tabularnewline
 &  & $W_{3}=2665230577$ & $\hgt(\frac{-18396006}{235225},\frac{300932687328}{114084125})$\tabularnewline
\hline 
$103$ & $(2,-3,-6)$ & $W_{0}=11581$ & $(-108,2772)$\tabularnewline
 &  & $W_{1}=462$ & $(10609/9,816893/27)$\tabularnewline
 &  & $W_{2}=13345$ & $(74778/49,16804656/343)$\tabularnewline
 &  & $W_{3}=1487$ & $(-10094/121,3564624/1331)$\tabularnewline
\hline 
$103$ & $(2,-6,-3)$ & $W_{0}=8487373$ & $(-103/841,-3044783/24389)$\tabularnewline
 &  & $W_{1}=16646$ & $(1039476,1059584484)$\tabularnewline
 &  & $W_{2}=8490735$ & $(1039682/1681,34458032/68921)$\tabularnewline
 &  & $W_{3}=8477279$ & $(-10086/49,57072/343)$\tabularnewline
\hline 
$198$ & $(1,2,2)$ & $W_{0}=58653195191109140161$ & $\hgt(\frac{28710309938}{15499969},\frac{3196495762009760}{61023377953})$\tabularnewline
 &  & $W_{1}=1573075476879053140$ & $\hgt(\frac{-3682885634214}{14355154969},\frac{12478518082656282720}{1719934182300797})$\tabularnewline
 &  & $W_{2}=66521235373358303439$ & $\hgt(\frac{-253808255431}{2179956100},\frac{667116748338909653}{101782150309000})$\tabularnewline
 &  & $W_{3}=22035538516500689039$ & $\hgt(\frac{5205735166800}{1275418369},\frac{10465997585174146680}{45549016212097})$\tabularnewline
\hline 
$198$ & $(2,-3,-6)$ & $W_{0}=255711950171342941$ & $\hgt(\frac{-496910700}{2436721},\frac{28292464294380}{3803721481})$\tabularnewline
 &  & $W_{1}=7360756127254530$ & $\hgt(\frac{96496588321}{41409225},\frac{22649052322875419}{266468362875})$\tabularnewline
 &  & $W_{2}=294877147817303041$ & $\hgt(\frac{677864759226}{236452129},\frac{454954755901005360}{3635924387633})$\tabularnewline
 &  & $W_{3}=26397138616197359$ & $\hgt(\frac{-94107947342}{567725929},\frac{97869913424403120}{13527205710283})$\tabularnewline
\hline 
$198$ & $(2,-6,-3)$ & $W_{0}=9901$ & $(-6,1680)$\tabularnewline
 &  & $W_{1}=70$ & $(79202,22176560)$\tabularnewline
 &  & $W_{2}=9999$ & $(59700/49,2376060/343)$\tabularnewline
 &  & $W_{3}=9601$ & $(-9751/25,277207/125)$\tabularnewline
\hline 
\end{tabular}

\end{document}